# RODEO: SPARSE, GREEDY NONPARAMETRIC REGRESSION

By John Lafferty[1] and Larry Wasserman[2]

*Carnegie Mellon University*

We present a greedy method for simultaneously performing local bandwidth selection and variable selection in nonparametric regression. The method starts with a local linear estimator with large bandwidths, and incrementally decreases the bandwidth of variables for which the gradient of the estimator with respect to bandwidth is large. The method—called *rodeo* (regularization of derivative expectation operator)—conducts a sequence of hypothesis tests to threshold derivatives, and is easy to implement. Under certain assumptions on the regression function and sampling density, it is shown that the rodeo applied to local linear smoothing avoids the curse of dimensionality, achieving near optimal minimax rates of convergence in the number of relevant variables, as if these variables were isolated in advance.

**1. Introduction.** Estimating a high-dimensional regression function is notoriously difficult due to the curse of dimensionality. Minimax theory precisely characterizes the curse. Let

$$Y_i = m(X_i) + \varepsilon_i, \qquad i = 1, \dots, n, \tag{1.1}$$

where $X_i = (X_i(1), \dots, X_i(d)) \in \mathbb{R}^d$ is a $d$-dimensional covariate, $m : \mathbb{R}^d \to \mathbb{R}$ is the unknown function to estimate and $\varepsilon_i \sim N(0, \sigma^2)$. Then if $m$ is in $W_2(c)$, the $d$-dimensional Sobolev ball of order two and radius $c$, it is well known that

$$\liminf_{n\to\infty} n^{4/(4+d)} \inf_{\widehat{m}_n} \sup_{m \in W_2(c)} \mathcal{R}(\widehat{m}_n, m) > 0, \tag{1.2}$$

Received June 2005; revised April 2007.

[1]Supported by NSF Grants CCR-01-22481, IIS-03-12814 and IIS-04-27206.

[2]Supported by NIH Grants R01-CA54852-07 and MH57881 and NSF Grant DMS-01-04016.

where $\mathcal{R}(\widehat{m}_n, m) = \mathbb{E}_m \int (\widehat{m}_n(x) - m(x))^2\, dx$ is the risk of the estimate $\widehat{m}_n$ constructed from a sample of size $n$ (Gyorfi et al. [12] and Stone et al. [25]). Thus, the best rate of convergence is $n^{-4/(4+d)}$, which is impractically slow if $d$ is large.

However, for some applications it is reasonable to expect that the true function only depends on a small number of the total covariates. Suppose that $m$ satisfies such a sparseness condition, so that

$$m(x) = m(x_R), \tag{1.3}$$

where $x_R = (x_j : j \in R)$, $R \subset \{1, \ldots, d\}$ is a subset of the $d$ covariates, of size $r = |R| \ll d$. We call $\{x_j\}_{j \in R}$ the *relevant variables*. Note that if an oracle were to identify and isolate the relevant variables, the better minimax rate of $n^{-4/(4+r)}$ could be achieved, and this would be the fastest rate possible. Thus, we are faced with the problem of variable selection in nonparametric regression. Our strategy is to seek a greedy method that incrementally searches through bandwidths in small steps.

A large body of previous work has addressed this fundamental problem, which has led to a variety of methods to combat the curse of dimensionality. Many of these are based on very clever, though often heuristic techniques. For additive models of the form $m(x) = \sum_j m_j(x_j)$, standard methods like stepwise selection, $C_p$ and AIC can be used (Hastie, Tibshirani and Friedman [14]). For spline models, Zhang et al. [31] use likelihood basis pursuit, essentially the lasso adapted to the spline setting. CART (Breiman et al. [1]) and MARS (Friedman [8]) effectively perform variable selection as part of their function fitting. Support vector regression can be seen as creating a sparse representation using basis pursuit in a reproducing kernel Hilbert space (Girosi [11]). There is also a large literature on Bayesian methods, including methods for sparse Gaussian processes (Tipping [27], Smola and Bartlett [24], Lawrence, Seeger and Herbrich [17]); see George and McCulloch [10] for a brief survey. More recently, Li, Cook and Nachsteim [19] use independence testing for variable selection and [2] introduced a boosting approach. While these methods have met with varying degrees of empirical success, they can be challenging to implement and demanding computationally. Moreover, these methods are typically very difficult to analyze theoretically, and so come with limited formal guarantees. Indeed, the theoretical analysis of sparse *parametric* estimators such as the lasso (Tibshirani [26]) is challenging, and only recently has significant progress been made on this front in the statistics and signal processing communities (Donoho [3], Fu and Knight [9], Tropp [28, 29], Fan and Peng [7] and Fan and Li [6]).

In this paper, we present a new approach for sparse nonparametric function estimation that is both computationally simple and amenable to theoretical analysis. We call the general framework *rodeo*, for "regularization of derivative expectation operator." It is based on the idea that bandwidth and variable selection can be simultaneously performed by computing the infinitesimal change in a nonparametric estimator as a function of the smoothing parameters, and then thresholding these derivatives to get a sparse estimate. As a simple version of this principle, we use hard thresholding, effectively carrying out a sequence of hypothesis tests. A modified version that replaces testing with soft thresholding may be viewed as solving a sequence of lasso problems. The potential appeal of this approach is that it can be based on relatively simple and theoretically well-understood nonparametric techniques such as local linear smoothing, leading to methods that are simple to implement and can be used in high-dimensional problems. Moreover, we show that they can achieve near optimal minimax rates of convergence, and therefore circumvent the curse of dimensionality when the true function is indeed sparse. When applied in one dimension, our method yields a local bandwidth selector and is similar to the estimators of Ruppert [21] and Lepski, Mammen and Spokoiny [18]. The method in Lepski, Mammen and Spokoiny [18] and its multivariate extension in Kerkyacharian, Lepski and Picard [16] yield estimators that are more refined than our method in the sense that their estimator is spatially adaptive over large classes of function spaces. However, their method is not greedy: it involves searching over a large class of bandwidths. Our goal is to develop a greedy method that scales to high dimensions.

Our method is related to the structural adaptation method of Hristache et al. [15] and Samarov, Spokoiny and Vial [23], which is designed for multi-index models. The general multi-index model is

$$Y = g_0(Tx) + \varepsilon, \tag{1.4}$$

where $x \in \mathbb{R}^d$ and $T$ is a linear orthonormal mapping from $\mathbb{R}^d$ onto $\mathbb{R}^r$ with $r < d$. Variable selection corresponds to taking $T$ to be a $r$ by $d$ matrix of 0's and 1's with each $T_{ij} = 1$ if $x_j$ is the $i$th relevant variable. Nonparametric variable selection can also be regarded as a special case of the partially linear model in Samarov, Spokoiny and Vial [23], which takes

$$Y = \theta^{\mathrm{T}} x_1 + G(x_2) + \varepsilon, \tag{1.5}$$

where $x = (x_1, x_2)$. Taking $\theta$ to be zero yields the model in this paper. The advantage of structural adaptation is that it yields, under certain conditions, $\sqrt{n}$ estimates of the image of $T$ in (1.4) and $\theta$ in (1.5). However, structural adaptation does not yield optimal bandwidths or optimal estimates of the regression function, although this is not the intended goal of the method.

In the following section we outline the basic rodeo approach, which is actually a general strategy that can be applied to a wide range of nonparametric estimators. We then specialize in Section 3 to the case of local linear smoothing, since the asymptotic properties of this smoothing technique are fairly well understood. In particular, we build upon the analysis of Ruppert and Wand [22] for local linear regression; a notable difference is that we allow the dimension to increase with sample size, which requires a more detailed analysis of the asymptotics. In Section 4 we present some simple examples of the rodeo, before proceeding to an analysis of its properties in Section 5. Our main theoretical result characterizes the asymptotic running time, selected bandwidths, and risk of the algorithm. Finally, in Section 6, we present further examples and discuss several extensions of the basic version of the rodeo considered in the earlier sections. The proofs of the theoretical properties of the rodeo are given in Section 7.

**2. Rodeo: the main idea.** The key idea in our approach is as follows. Fix a point $x$ and let $\widehat{m}_h(x)$ denote an estimator of $m(x)$ based on a vector of smoothing parameters $h = (h_1, \ldots, h_d)$. If $c$ is a scalar, then we write $h = c$ to mean $h = (c, \ldots, c)$.

Let $M(h) = \mathbb{E}(\widehat{m}_h(x))$ denote the mean of $\widehat{m}_h(x)$. For now, assume that $x = x_i$ is one of the observed data points and that $\widehat{m}_0(x) = Y_i$. In that case, $m(x) = M(0) = \mathbb{E}(Y_i)$. If $P = (h(t) : 0 \leq t \leq 1)$ is a smooth path through the set of smoothing parameters with $h(0) = 0$ and $h(1) = 1$ (or any other fixed, large bandwidth) then

$$m(x) = M(0) = M(1) + M(0) - M(1) \tag{2.1a}$$

$$= M(1) - \int_0^1 \frac{dM(h(s))}{ds}\, ds \tag{2.1b}$$

$$= M(1) - \int_0^1 \langle D(h(s)), \dot{h}(s)\rangle\, ds, \tag{2.1c}$$

where

$$D(h) = \nabla M(h) = \left(\frac{\partial M}{\partial h_1}, \ldots, \frac{\partial M}{\partial h_d}\right)^{\mathrm{T}} \tag{2.2}$$

is the gradient of $M(h)$ and $\dot{h}(s) = \frac{dh(s)}{ds}$ is the derivative of $h(s)$ along the path. An unbiased, low variance estimator of $M(1)$ is $\widehat{m}_1(x)$. An unbiased estimator of $D(h)$ is

$$Z(h) = \left(\frac{\partial \widehat{m}_h(x)}{\partial h_1}, \ldots, \frac{\partial \widehat{m}_h(x)}{\partial h_d}\right)^{\mathrm{T}}. \tag{2.3}$$

The naive estimator

$$\widehat{m}(x) = \widehat{m}_1(x) - \int_0^1 \langle Z(h(s)), \dot{h}(s)\rangle\, ds \tag{2.4}$$

is identically equal to $\widehat{m}_0(x) = Y_i$, which has poor risk since the variance of $Z(h)$ is large for small $h$. However, our sparsity assumption on $m$ suggests that there should be paths for which $D(h)$ is also sparse. Along such a path, we replace $Z(h)$ with an estimator $\widehat{D}(h)$ that makes use of the sparsity assumption. Our estimate of $m(x)$ is then

$$\widetilde{m}(x) = \widehat{m}_1(x) - \int_0^1 \langle \widehat{D}(s), \dot{h}(s) \rangle \, ds. \tag{2.5}$$

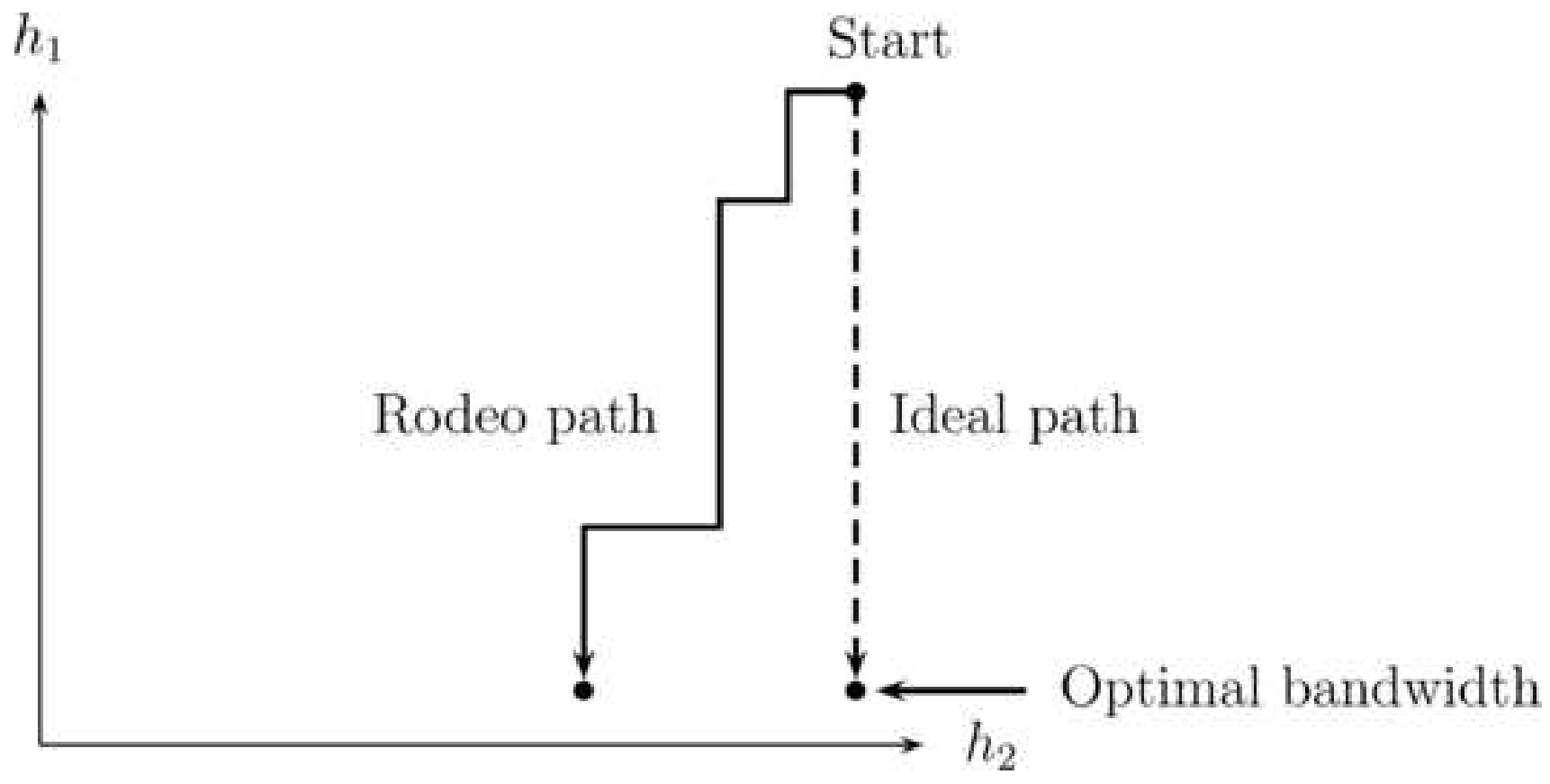

FIG. 1. *Conceptual illustration: The bandwidths for the relevant variables ($h_1$) are shrunk, while the bandwidths for the irrelevant variables ($h_2$) are kept relatively large.*

To implement this idea we need to do two things: (i) we need to find a path for which the derivative is sparse and (ii) we need to take advantage of this sparseness when estimating $D$ along that path.

The key observation is that if $x_j$ is irrelevant, then we expect that changing the bandwidth $h_j$ for that variable should cause only a small change in the estimator $\widehat{m}_h(x)$. Conversely, if $x_j$ is relevant, then we expect that changing the bandwidth $h_j$ for that variable should cause a large change in the estimator. Thus, $Z_j(h) = \partial \widehat{m}_h(x)/\partial h_j$ should discriminate between relevant and irrelevant covariates. To simplify the procedure, we can replace the continuum of bandwidths in the interval with a discrete set where each $h_j \in \mathcal{B} = \{h_0, \beta h_0, \beta^2 h_0, \dots\}$ for some $0 < \beta < 1$. Moreover, we can proceed in a greedy fashion by estimating $D(h)$ sequentially with $h_j \in \mathcal{B}$ and setting $\widehat{D}_j(h) = 0$ when $h_j < \widehat{h}_j$, where $\widehat{h}_j$ is the first $h$ such that $|Z_j(h)| < \lambda_j(h)$ for some threshold $\lambda_j$. This greedy version, coupled with the hard threshold estimator, yields $\widetilde{m}(x) = \widehat{m}_{\widehat{h}}(x)$. A conceptual illustration of the idea is shown in Figure 1.

To further elucidate the idea, consider now the one-dimensional case $x \in \mathbb{R}$, so that

$$m(x) = M(1) - \int_0^1 \frac{dM(h)}{dh} \, dh = M(1) - \int_0^1 D(h) \, dh. \tag{2.6}$$

Suppose that $\widehat{m}_h(x) = \sum_{i=1}^n Y_i \ell_i(x,h)$ is a linear estimator, where the weights $\ell_i(x,h)$ depend on a bandwidth $h$.

In this case

$$Z(h) = \sum_{i=1}^n Y_i \ell_i'(x,h) \tag{2.7}$$

where the prime denotes differentiation with respect to $h$. Then we set

$$\widetilde{m}(x) = \widehat{m}_1(x) - \int_0^1 \widehat{D}(h)\,dh \tag{2.8}$$

where $\widehat{D}(h)$ is an estimator of $D(h)$. Now,

$$Z(h) \approx N(b(h), s^2(h)) \tag{2.9}$$

where, for typical smoothers, $b(h) \approx Ah$ and $s^2(h) \approx C/nh^3$ for some constants $A$ and $C$. Take the hard threshold estimator

$$\widehat{D}(h) = Z(h) I(|Z(h)| > \lambda(h)), \tag{2.10}$$

where $\lambda(h)$ is chosen to be slightly larger than $s(h)$. An alternative is the soft-threshold estimator

$$\widehat{D}(h) = \operatorname{sign}(Z(h))(|Z(h)| - \lambda(h))_+. \tag{2.11}$$

The greedy algorithm, coupled with the hard threshold estimator, yields a bandwidth selection procedure based on testing. This approach to bandwidth selection is very similar to that of Lepski, Mammen and Spokoiny [18], who take

$$\widehat{h} = \max\{h \in \mathcal{H} : \phi(h,\eta) = 0 \text{ for all } \eta < h\}, \tag{2.12}$$

where $\phi(h,\eta)$ is a test for whether $\widehat{m}_\eta$ improves on $\widehat{m}_h$. This more refined test leads to estimators that achieve good spatial adaptation over large function classes. Kerkyacharian, Lepski and Picard [16] extend the idea to multiple dimensions. Our approach is also similar to a method of Ruppert [21] that uses a sequence of decreasing bandwidths and then estimates the optimal bandwidth by estimating the mean squared error as a function of bandwidth. Our greedy approach only tests whether an infinitesimal change in the bandwidth from its current setting leads to a significant change in the estimate, and is more easily extended to a practical method in higher dimensions.

**3. Rodeo using local linear regression.** Now we present the multivariate rodeo in detail. We use local linear smoothing as the basic method since it is known to have many good properties. Let $x = (x(1), \ldots, x(d))$ be some target point at which we want to estimate $m$. Let $\widehat{m}_H(x)$ denote the local linear estimator of $m(x)$ using bandwidth matrix $H$. Thus,

$$\widehat{m}_H(x) = e_1^{\mathrm{T}} (X_x^{\mathrm{T}} W_x X_x)^{-1} X_x^{\mathrm{T}} W_x Y \equiv S_x Y, \tag{3.1}$$

where $e_1 = (1, 0, \ldots, 0)^{\mathrm{T}}$,

$$X_x = \begin{pmatrix} 1 & (X_1 - x)^{\mathrm{T}} \\ \vdots & \vdots \\ 1 & (X_n - x)^{\mathrm{T}} \end{pmatrix}, \tag{3.2}$$

$W_x$ is diagonal with $(i, i)$ element $K_H(X_i - x)$ and $K_H(u) = |H|^{-1/2} K(H^{-1/2}u)$. The estimator $\widehat{m}_H$ can be written as

$$\widehat{m}_H(x) = \sum_{i=1}^{n} G(X_i, x, h) Y_i, \tag{3.3}$$

where

$$G(u, x, h) = e_1^{\mathrm{T}} (X_x^{\mathrm{T}} W_x X_x)^{-1} \begin{pmatrix} 1 \\ (u - x)^{\mathrm{T}} \end{pmatrix} K_H(u - x) \tag{3.4}$$

is called the *effective kernel*. One can regard local linear regression as a refinement of kernel regression where the effective kernel $G$ adjusts for boundary bias and design bias; see Fan [5], Hastie and Loader [13] and Ruppert and Wand [22].

We assume that the covariates are random with density $f(x)$ and that $x$ is interior to the support of $f$. We make the same assumptions as Ruppert and Wand [22] in their analysis of the bias and variance of local linear regression. In particular:

(i) The kernel $K$ has compact support with zero odd moments and there exists $\nu_2 = \nu_2(K) \neq 0$ such that

$$\int u u^{\mathrm{T}} K(u)\, du = \nu_2(K) I, \tag{3.5}$$

where $I$ is the $d \times d$ identity matrix.

(ii) The sampling density $f(x)$ is continuously differentiable and strictly positive.

In the version of the algorithm that follows, we take $K$ to be a product kernel and $H$ to be diagonal with elements $h = (h_1, \ldots, h_d)$ and we write $\widehat{m}_h$ instead of $\widehat{m}_H$.

Our method is based on the statistic

$$Z_j = \frac{\partial \widehat{m}_h(x)}{\partial h_j} = \sum_{i=1}^{n} G_j(X_i, x, h) Y_i, \tag{3.6}$$

where

$$G_j(u, x, h) = \frac{\partial G(u, x, h)}{\partial h_j}. \tag{3.7}$$

*Rodeo: Hard thresholding version*

1. *Select* constant $0 < \beta < 1$ and initial bandwidth

$$h_0 = \frac{c_0}{\log \log n}. \tag{3.11}$$

2. *Initialize* the bandwidths, and activate all covariates:
   (a) $h_j = h_0$, $j = 1, 2, \ldots, d$.
   (b) $\mathcal{A} = \{1, 2, \ldots, d\}$.

3. *While $\mathcal{A}$ is nonempty*, do for each $j \in \mathcal{A}$:
   (a) Compute the estimated derivative expectation: $Z_j$ [equation (3.6)] and $s_j$ [equation (3.9)].
   (b) Compute the threshold $\lambda_j = s_j \sqrt{2 \log n}$.
   (c) If $|Z_j| > \lambda_j$, then set $h_j \leftarrow \beta h_j$; otherwise remove $j$ from $\mathcal{A}$.

4. *Output* bandwidths $h^\star = (h_1, \ldots, h_d)$ and estimator $\widetilde{m}(x) = \widehat{m}_{h^\star}(x)$.

FIG. 2. *The hard thresholding version of the rodeo, which can be applied using the derivatives $Z_j$ of any nonparametric smoother. The algorithm stops when all derivatives are below threshold. As shown in the theoretical analysis, this happens after $T_n = O(\log n)$ iterations.*

Let

$$\mu_j \equiv \mu_j(h) = \mathbb{E}(Z_j | X_1, \ldots, X_n) = \sum_{i=1}^n G_j(X_i, x, h) m(X_i) \tag{3.8}$$

and

$$\operatorname{Var}(Z_j | X_1, \ldots, X_n) = \sigma^2 \sum_{i=1}^n G_j(X_i, x, h)^2. \tag{3.9}$$

In Section 4.3 we explain how to estimate $\sigma$; for now, assume that $\sigma$ is known. The hard thresholding version of the rodeo algorithm is described in Figure 2.

To derive an explicit expression for $Z_j$, equivalently $G_j$, we use

$$\frac{\partial A^{-1}}{\partial h} = -A^{-1} \frac{\partial A}{\partial h} A^{-1} \tag{3.10}$$

to get that

$$\begin{aligned} Z_j &= \frac{\partial \widehat{m}_h(x)}{\partial h_j} \\ &= e_1^{\mathrm{T}} (X^{\mathrm{T}} W X)^{-1} X^{\mathrm{T}} \frac{\partial W}{\partial h_j} Y \end{aligned} \tag{3.12a}$$

$$- e_1^{\mathrm{T}}(X^{\mathrm{T}}WX)^{-1}X^{\mathrm{T}}\frac{\partial W}{\partial h_j}X(X^{\mathrm{T}}WX)^{-1}X^{\mathrm{T}}WY \tag{3.12b}$$

$$= e_1^{\mathrm{T}}(X^{\mathrm{T}}WX)^{-1}X^{\mathrm{T}}\frac{\partial W}{\partial h_j}(Y - X\widehat{\alpha}), \tag{3.12c}$$

where $\widehat{\alpha} = (X^{\mathrm{T}}WX)^{-1}X^{\mathrm{T}}WY$ is the coefficient vector for the local linear fit (and we have dropped the dependence on the local point $x$ in the notation).

Note that the factor $|H|^{-1} = \prod_{i=1}^{d} 1/h_i$ in the kernel cancels in the expression for $\widehat{m}$, and therefore we can ignore it in our calculation of $Z_j$. Assuming a product kernel we have

$$W = \mathrm{diag}\left(\prod_{j=1}^{d} K((X_{1j} - x_j)/h_j), \ldots, \prod_{j=1}^{d} K((X_{nj} - x_j)/h_j)\right) \tag{3.13}$$

and $\partial W/\partial h_j = WL_j$, where

$$L_j = \mathrm{diag}\left(\frac{\partial \log K((X_{1j} - x_j)/h_j)}{\partial h_j}, \ldots, \frac{\partial \log K((X_{nj} - x_j)/h_j)}{\partial h_j}\right), \tag{3.14}$$

and thus

$$\begin{aligned} Z_j &= e_1^{\mathrm{T}}(X^{\mathrm{T}}WX)^{-1}X^{\mathrm{T}}WL_j(Y - X\widehat{\alpha}) \\ &= e_1^{\mathrm{T}}BL_j(I - XB)Y = G_j(x,h)^{\mathrm{T}}Y \end{aligned} \tag{3.15}$$

where $B = (X^{\mathrm{T}}WX)^{-1}X^{\mathrm{T}}W$.

The calculation of $L_j$ is typically straightforward. As two examples, with the Gaussian kernel $K(u) = \exp(-u^2/2)$ we have

$$L_j = \frac{1}{h_j^3}\mathrm{diag}((X_{1j} - x_j)^2, \ldots, (X_{nj} - x_j)^2) \tag{3.16}$$

and for the Epanechnikov kernel $K(u) = (5 - x^2)\mathbb{I}(|x| \le \sqrt{5})$ we have

$$L_j = \frac{1}{h_j^3}\mathrm{diag}\left(\frac{2(X_{1j} - x_j)^2}{5 - (X_{1j} - x_j)^2/h_j^2}\mathbb{I}(|X_{1j} - x_j| \le \sqrt{5}h_j), \ldots,\right. \tag{3.17a}$$

$$\left.\frac{2(X_{nj} - x_j)^2}{5 - (X_{nj} - x_j)^2/h_j^2}\mathbb{I}(|X_{1j} - x_j| \le \sqrt{5}h_j)\right). \tag{3.17b}$$

**4. Examples.** In this section we illustrate the rodeo on some examples. We return to the examples later when we discuss estimating $\sigma$, as well as a global (nonlocal) version of the rodeo.

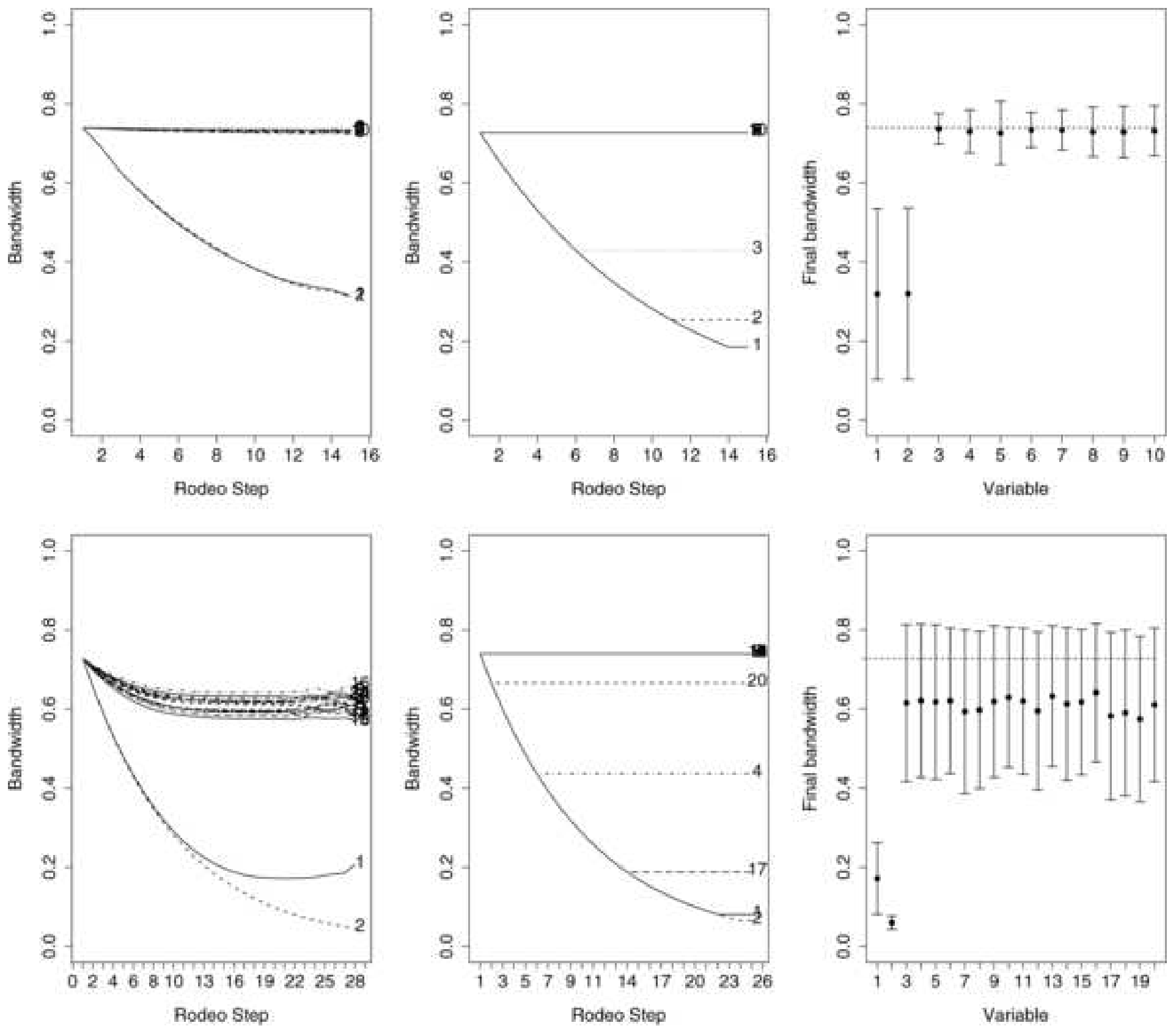

Fig. 3. *Rodeo run on synthetic data sets, showing average bandwidths over 200 runs (left), final bandwidths with standard errors (right), and bandwidths on a single run of the algorithm (center). In the top plots the regression function is* $m(x) = 5x_1^2x_2^2$ *with* $d = 10$, $n = 500$ *and* $\sigma = 0.5$ *and in the lower plots the regression function is* $m(x) = 2(x_1 + 1)^3 + 2\sin(10x_2)$, $d = 20$, $n = 750$ *and* $\sigma = 1$. *The figures show that the bandwidths for the relevant variables* $x_1$ *and* $x_2$ *are shrunk, while the bandwidths for the irrelevant variables remain large.*

4.1. *Two relevant variables.* In the first example, we take $m(x) = 5x_1^2x_2^2$ with $d = 10$, $\sigma = 0.5$ with $x_i \sim \text{Uniform}(0, 1)$. The algorithm is applied to the local linear estimates around the test point $x_0 = (\frac{1}{2}, \ldots, \frac{1}{2})$, with $\beta = 0.8$. Figure 3 shows the bandwidths averaged over 200 runs of the rodeo, on data sets of size $n = 750$. The second example in Figure 4 shows the algorithm applied to the function $m(x) = 2(x_1 + 1)^3 + 2\sin(10x_2)$, in this case in $d = 20$ dimensions with $\sigma = 1$.

The plots demonstrate how the bandwidths $h_1$ and $h_2$ of the relevant variables are shrunk, while the bandwidths of the relevant variables tend to remain large.

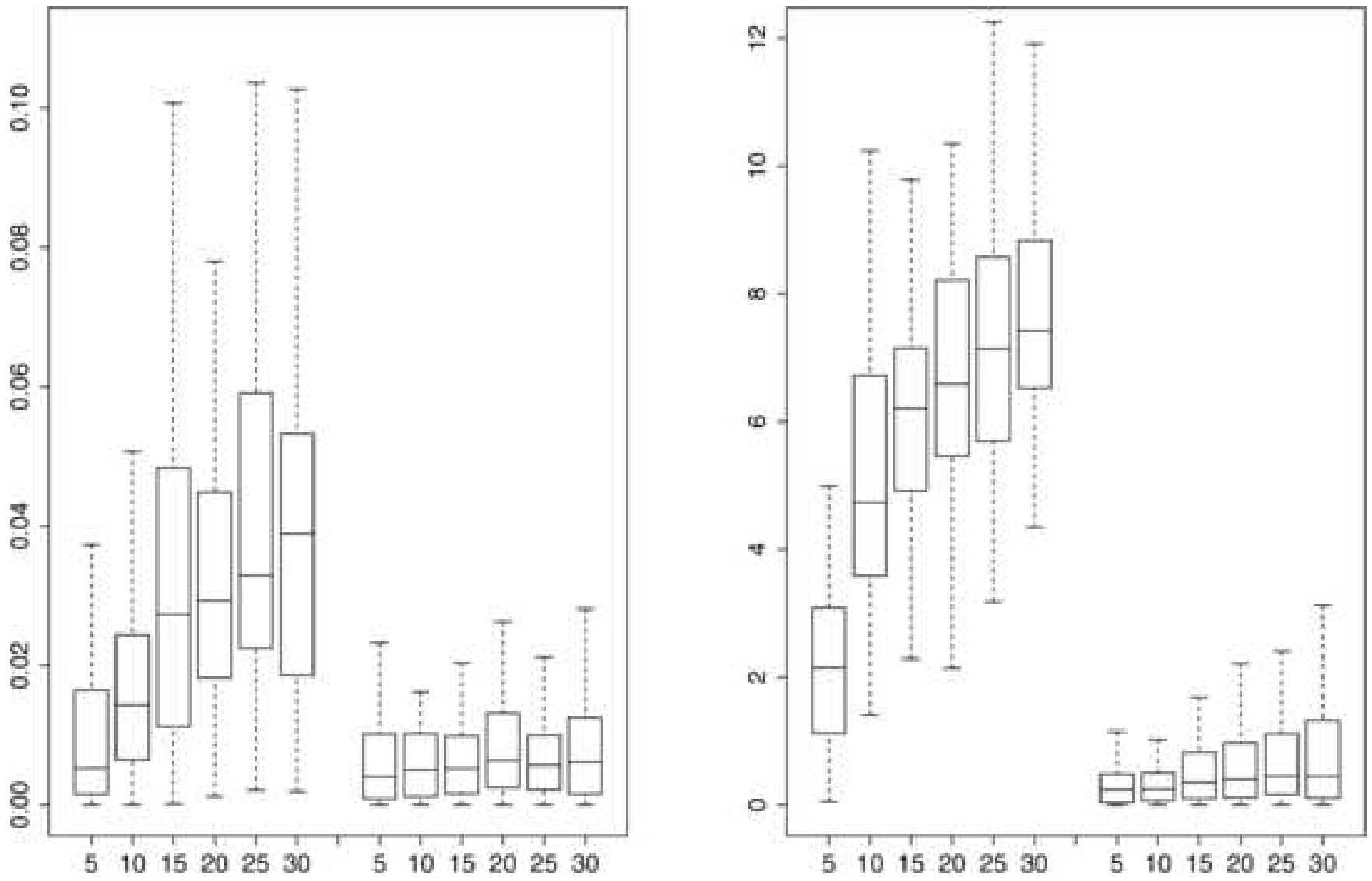

FIG. 4. *Squared error of the estimator on the previous examples,* $m(x) = 5x_1^2x_2^2$ *(left) and* $m(x) = 2(x_1+1)^3 + 2\sin(10x_2)$ *(right). For each plot, the left six boxplots show the risk in different dimensions* ($d = 5, 10, 15, 20, 25, 30$) *when using a single bandwidth, chosen by leave-one-out cross validation. The right six boxplots show the squared error on the same data with bandwidths selected using the rodeo.*

4.2. *A one-dimensional example.* Figure 5 illustrates the algorithm in one dimension. The underlying function in this case is $m(x) = (1/x)\sin(15/x)$, and $n = 1{,}500$ data points are sampled as $x \sim \text{Uniform}(0,1) + \frac{1}{2}$. The algorithm is run at two test points; the function is more rapidly varying near the test point $x = 0.67$ than near the test point $x = 1.3$, and the rodeo appropriately selects a smaller bandwidth at $x = 0.67$. The right plot of Figure 5 displays boxplots for the logarithm of the final bandwidth in the base $1/\beta$ (equivalently, minus the number of steps in the algorithm) where $\beta = 0.8$, averaged over 50 randomly generated data sets.

The figure illustrates how smaller bandwidths are selected where the function is more rapidly varying. However, we do not claim that the method is adaptive over large classes of function spaces. As discussed earlier, the technique is intentionally a greedy algorithm; adapting to unknown smoothness may require a more refined search over bandwidths that does not scale to large dimensions, and is out of the scope of the current paper.

4.3. *Estimating* $\sigma$. The algorithm requires that we insert an estimate $\widehat{\sigma}$ of $\sigma$ in (3.9). An estimator for $\sigma$ can be obtained by generalizing a method of Rice [20]. For $i < \ell$, let

$$d_{i\ell} = \|X_i - X_\ell\|. \tag{4.1}$$

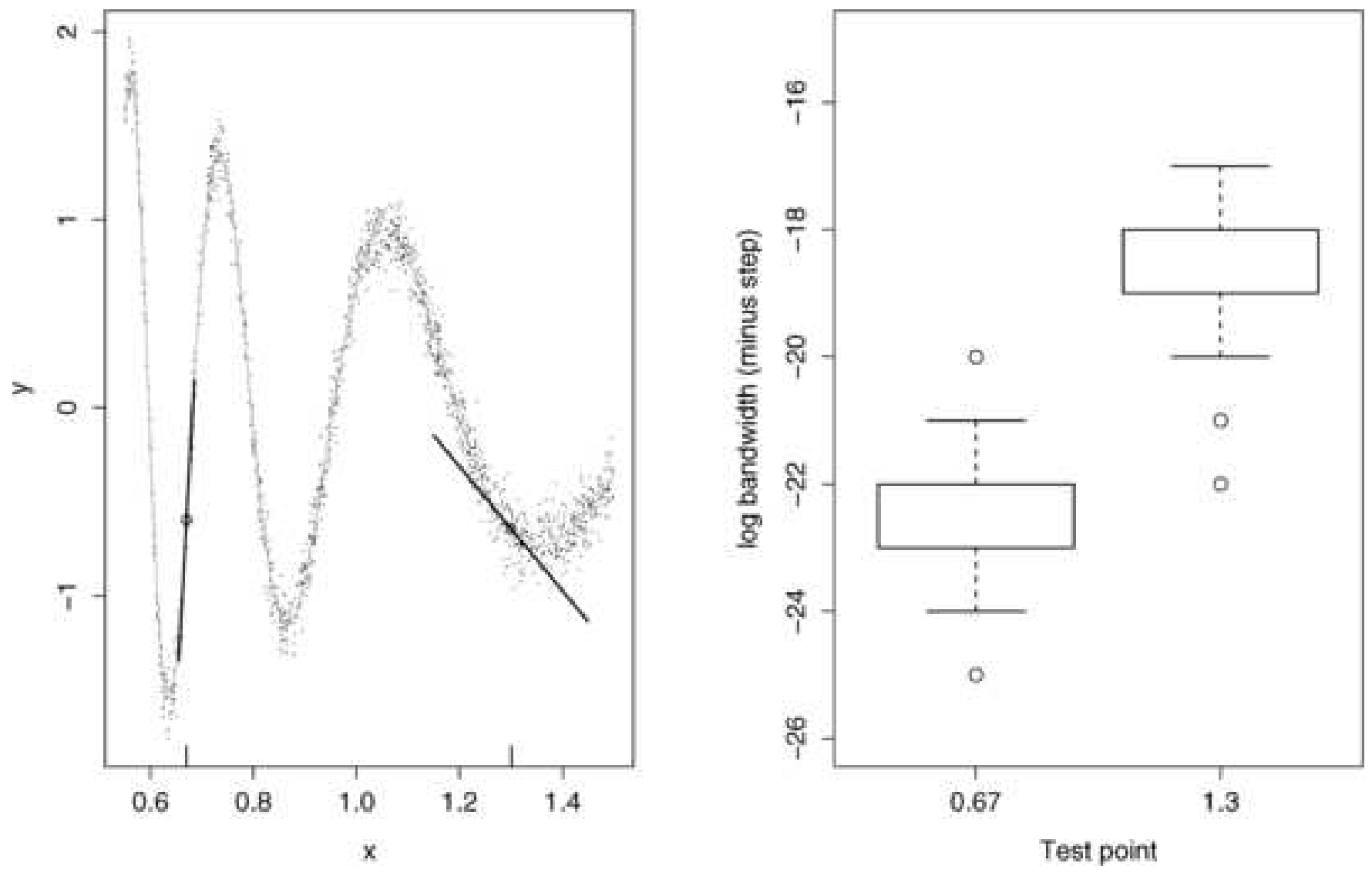

FIG. 5. *A one-dimensional example. The regression function is* $m(x) = (1/x)\sin(15/x)$, *and* $n = 1{,}500$ *data points are sampled,* $x \sim \text{Uniform}(0,1) + \frac{1}{2}$. *The left plot shows the local linear fit at two test points; the right plot shows the final log bandwidth,* $\log_{1/\beta} h_\star$ *(equivalently, minus the number of steps) of the rodeo over 50 randomly generated data sets.*

Fix an integer $J$ and let $\mathcal{E}$ denote the set of pairs $(i,\ell)$ corresponding the $J$ smallest values of $d_{i\ell}$. Now define

$$\widehat{\sigma}^2 = \frac{1}{2J} \sum_{i,\ell \in \mathcal{E}} (Y_i - Y_\ell)^2. \tag{4.2}$$

Then

$$\mathbb{E}(\widehat{\sigma}^2) = \sigma^2 + \text{bias}, \tag{4.3}$$

where

$$\text{bias} \le D \sup_x \sum_{j=1}^{r} \left| \frac{\partial m_j(x)}{\partial x_j} \right| \tag{4.4}$$

with $D$ given by

$$D = \max_{i,\ell \in \mathcal{E}} \|X_i - X_\ell\|. \tag{4.5}$$

There is a bias-variance tradeoff: large $J$ makes $\widehat{\sigma}^2$ positively biased, and small $J$ makes $\widehat{\sigma}^2$ highly variable. Note, however, that the bias is mitigated by sparsity (small $r$).

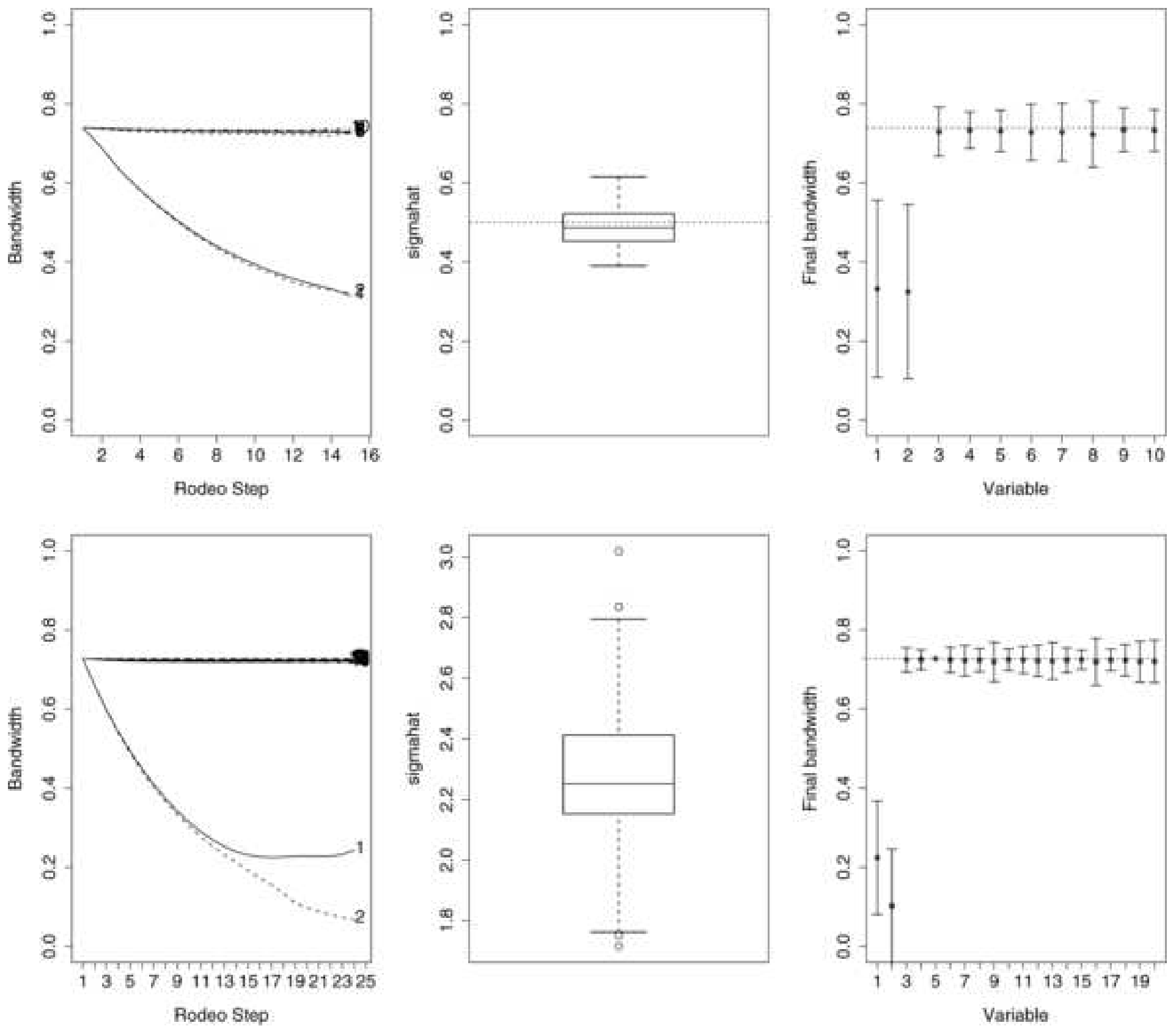

Fig. 6. *Rodeo run on the examples of Section 4.1, but now estimating the noise using the estimate $\widehat{\sigma}$ discussed in Section 4.3. Top: $\sigma = 0.5$, $d = 10$; bottom: $\sigma = 1$, $d = 20$. In higher dimensions the noise is over-estimated (center plots), which results in the irrelevant variables being more aggressively eliminated; compare Figure 3.*

A more robust estimate may result from taking

$$\widehat{\sigma} = \frac{\sqrt{\pi}}{2} \operatorname{median}\{|Y_i - Y_\ell|\}_{i,\ell \in \mathcal{E}} \tag{4.6}$$

where the constant comes from observing that if $X_i$ is close to $X_\ell$, then

$$|Y_i - Y_\ell| \sim |N(0, 2\sigma^2)| = \sqrt{2}\sigma|Z|, \tag{4.7}$$

where $Z$ is a standard normal with $\mathbb{E}|Z| = \sqrt{2/\pi}$.

Now we redo the earlier examples, taking $\sigma$ as unknown. Figure 6 shows the result of running the algorithm on the examples of Section 4.1, however, now estimating the noise using estimate (4.6). For the higher-dimensional example, with $d = 20$, the noise variance is over-estimated, with the primary result that the irrelevant variables are more aggressively thresholded out; compare Figure 6 to Figure 3.

Although we do not pursue it in this paper, there is also the possibility of allowing $\sigma(x)$ to be a function of $x$ and estimating it locally.

4.4. *Computational cost.* When based on a local linear estimator, each step of the rodeo algorithm has the same computational cost as constructing a single local linear fit. This is dominated by the cost of constructing the matrix inverse $(X^{\mathrm{T}}WX)^{-1}$ in equation (3.15). Since the derivative needs to be computed for every variable, the algorithm thus scales as $O(d^4)$ in the dimension $d$. Implemented in R, the 20 dimensional example in Figure 3 takes 4 hours, 4 minutes and 40 seconds for 200 runs, or 73.4 seconds per run, when executed on a 1.5 GHz PowerPC Macintosh laptop. Although we focus on local linear regression, it should be noted that very similar results are obtained with kernel regression, which requires no matrix inversion. Using kernel regression, the same example requires 12 minutes and 33 seconds, or 3.7 seconds per run.

**5. Properties of the rodeo.** We now give some results on the statistical properties of the hard thresholding version of the rodeo estimator. Formally, we use a triangular array approach so that $m(x)$, $f(x)$, $d$ and $r$ can all change as $n$ changes, although we often suppress the dependence on $n$ for notational clarity. We assume throughout that $m$ has continuous third order derivatives in a neighborhood of $x$. For convenience of notation, we assume that the covariates are numbered such that the relevant variables $x_j$ correspond to $1 \leq j \leq r$ and the irrelevant variables $x_j$ correspond to $r+1 \leq j \leq d$.

A key aspect of our analysis is that we allow the dimension $d$ to increase with sample size $n$, and show that the algorithm achieves near optimal minimax rates of convergence if $d = O(\log n / \log\log n)$. This hinges on a careful analysis of the asymptotic bias and variance of the estimated derivative $Z_j$, taking the increasing dimension into account. We conjecture that, without further assumptions, $d$ cannot increase at a significantly faster rate, while obtaining near optimal rates of convergence.

The results are stated below, with the complete proofs given in Section 7.

Our main theoretical result characterizes the asymptotic running time, selected bandwidths and risk of the algorithm. In order to get a practical algorithm, we need to make assumptions on the functions $m$ and $f$.

(A1) The density $f(x)$ of $(X_1, \ldots, X_d)$ is uniform on the unit cube.

(A2)

$$\liminf_{n\to\infty} \min_{1\leq j\leq r} |m_{jj}(x)| > 0. \tag{5.1}$$

(A3) All derivatives of $m$ up to and including fourth order are bounded.

Assumption (A1) greatly simplifies the proofs. If we drop (A1), it is necessary to use a smaller starting bandwidth.

THEOREM 5.1. *Assume the conditions of Lemma 7.1 and suppose that assumptions* (A1) *and* (A2) *hold. In addition, suppose that*

$$A_{\min} = \min_{j \le r} |m_{jj}(x)| = \widetilde{\Omega}(1)$$

*and*

$$A_{\max} = \max_{j \le r} |m_{jj}(x)| = \widetilde{O}(1).$$

*Then the rodeo outputs bandwidths $h^\star$ that satisfy*

$$\mathbb{P}(h_j^\star = h_0 \text{ for all } j > r) \to 1 \tag{5.2}$$

*and for every $\varepsilon > 0$,*

$$\mathbb{P}(n^{-1/(4+r)-\varepsilon} \le h_j^\star \le n^{-1/(4+r)+\varepsilon} \text{ for all } j \le r) \to 1. \tag{5.3}$$

*Let $T_n$ be the stopping time of the algorithm. Then $\mathbb{P}(t_L \le T_n \le t_U) \to 1$ where*

$$t_L = \frac{1}{(r+4)\log(1/\beta)} \log\left(\frac{nA_{\min}^2}{8C^2 \log n (\log\log n)^d}\right), \tag{5.4}$$

$$t_U = \frac{1}{(r+4)\log(1/\beta)} \log\left(\frac{nA_{\max}^2}{aC^2 \log n (\log\log n)^d}\right) \tag{5.5}$$

*and $0 < a < 2$.*

COROLLARY 5.2. *Under the conditions of Theorem 5.1,*

$$(\widehat{m}_{h^\star}(x) - m(x))^2 = O_P(n^{-4/(4+r)+\varepsilon}) \tag{5.6}$$

*for every $\varepsilon > 0$.*

**6. Extensions and variations of the rodeo.** The rodeo represents a general strategy for nonparametric estimation, based on the idea of regularizing or testing the derivatives of an estimator with respect to smoothing parameters. There are many ways in which this basic strategy can be realized. In this section we discuss several variants of the basic hard thresholding version of the rodeo, including a soft thresholding version, a global rather than local bandwidth selection procedure, the use of testing and generalized cross validation, and connections to least angle regression. Further numerical examples are also given to illustrate these ideas.

*Rodeo: Soft thresholding version*

1. *Select* parameter $0 < \beta < 1$ and initial bandwidth $h_0$.
2. *Initialize* the bandwidths, and activate all covariates:
   (a) $h_j = h_0$, $j = 1, 2, \dots, d$.
   (b) $\mathcal{A} = \{1, 2, \dots, d\}$.
   (c) Initialize step, $t = 1$.
3. *While $\mathcal{A}$ is nonempty*
   (a) Set $dh_j(t) = 0$, $j = 1, \dots, d$.
   (b) Do for each $j \in \mathcal{A}$:
      (1) Compute the estimated derivative expectation $Z_j$ and $s_j$.
      (2) Compute the threshold $\lambda_j = s_j \sqrt{2 \log n}$.
      (3) If $|Z_j| > \lambda_j$, set $dh_j(t) = (1 - \beta) h_j$ and $h_j \leftarrow \beta h_j$; otherwise remove $j$ from $\mathcal{A}$.
      (4) Set $\widehat{D}_j(t) = \operatorname{sign}(Z_j(h))(|Z_j(h)| - \lambda_j)_+$.
   (c) Increment step, $t \leftarrow t + 1$.
4. *Output* bandwidths $h^\star = (h_1, \dots, h_d)$ and estimator

$$\widetilde{m}(x) = \widehat{m}_{h_0}(x) - \sum_{s=1}^{t} \langle \widehat{D}(s), dh(s) \rangle \tag{6.4}$$

FIG. 7. *The soft thresholding version of the rodeo.*

6.1. *Subtracting off a linear lasso.* Local linear regression is a nonparametric method that contains linear regression as a special case when $h \to \infty$. If the true function is linear but only a subset of the variables are relevant, then the rodeo will fail to separate the relevant and irrelevant variables since relevance is defined in terms of departures from the limiting parametric model. Indeed, the results depend on the Hessian of $m$ which is zero in the linear case. The rodeo may return a full linear fit with all variables. A simple modification can potentially fix this problem. First, do linear variable selection using, say, the lasso (Tibshirani [26]). Then run the rodeo on the residuals from that fit, but using all of the variables. An example of this procedure is given below in Section 6.4.

6.2. *Other estimators and other paths.* We have taken the estimate

$$\widehat{D}_j(h) = Z_j(h) I(|Z_j(h)| > \lambda_j) \tag{6.1}$$

with the result that

$$\widetilde{m}(x) = \widehat{m}_{h_0}(x) - \int_0^1 \langle \widehat{D}(s), \dot{h}(s) \rangle \, ds = \widehat{m}_{h^\star}(x). \tag{6.2}$$

There are many possible generalizations. First, we can replace $\widehat{D}$ with the soft-thresholded estimate

$$\widehat{D}_j(t) = \operatorname{sign}(Z_j(h))(|Z_j(h)| - \lambda_j)_+ \tag{6.3}$$

where the index $t$ denotes the $t$th step of the algorithm. Since $h_j$ is updated multiplicatively as $h_j \leftarrow \beta h_j$, the differential $dh_j(t)$ is given by $dh_j(t) = (1-\beta)h_j$. Using the resulting estimate of $D(t)$ and finite difference approximation for $\dot{h}(t)$ leads to the algorithm detailed in Figure 7.

Figure 8 shows a comparison of the hard and soft thresholding versions of the rodeo on the example function $m(x) = 2(x_1+1)^3 + 2\sin(10x_2)$ in $d = 10$ dimensions with $\sigma = 1$; $\beta$ was set to 0.9. For each of 100 randomly generated datasets, a random test point $x \sim \text{Uniform}(0,1)^d$ was generated, and the difference in losses was computed:

$$(\widetilde{m}_{\text{hard}}(x) - m(x))^2 - (\widetilde{m}_{\text{soft}}(x) - m(x))^2. \tag{6.5}$$

Thus, positive values indicate an advantage for soft thresholding, which is seen to be slightly more robust on this example.

Another natural extension would be to consider more general paths than paths that are restricted to be parallel to the axes. We leave this direction to future work.

6.3. *Global rodeo.* We have focused on estimation of $m$ locally at a point $x$. The idea can be extended to carry out global bandwidth and variable selection by averaging over multiple evaluation points $x_1,\dots,x_k$. These could be points of interest for estimation, could be randomly chosen, or could be taken to be identical to the observed $X_i$'s.

Averaging the $Z_j$'s directly leads to a statistic whose mean for relevant variables is asymptotically $k^{-1}h_j\sum_{i=1}^k m_{jj}(x_i)$. Because of sign changes in $m_{jj}(x)$, cancellations can occur resulting in a small value of the statistic for relevant variables. To eliminate the sign cancellation, we square the statistic. Another way of deriving a global method would be to use the statistic $\sup_x |Z_j^*(x)|$.

Let $x_1,\dots,x_k$ denote the evaluation points. Let

$$Z_j(x_i) = \sum_{s=1}^n Y_s G_j(X_s, x_i). \tag{6.6}$$

Then define the statistic

$$T_j \equiv \frac{1}{k}\sum_{i=1}^k Z_j^2(x_i) = \frac{1}{k} Y^{\mathrm{T}} P_j Y, \tag{6.7}$$

where $P_j = \mathcal{G}_j\mathcal{G}_j^{\mathrm{T}}$, with $\mathcal{G}_j(s,i) = G_j(X_s, x_i)$.

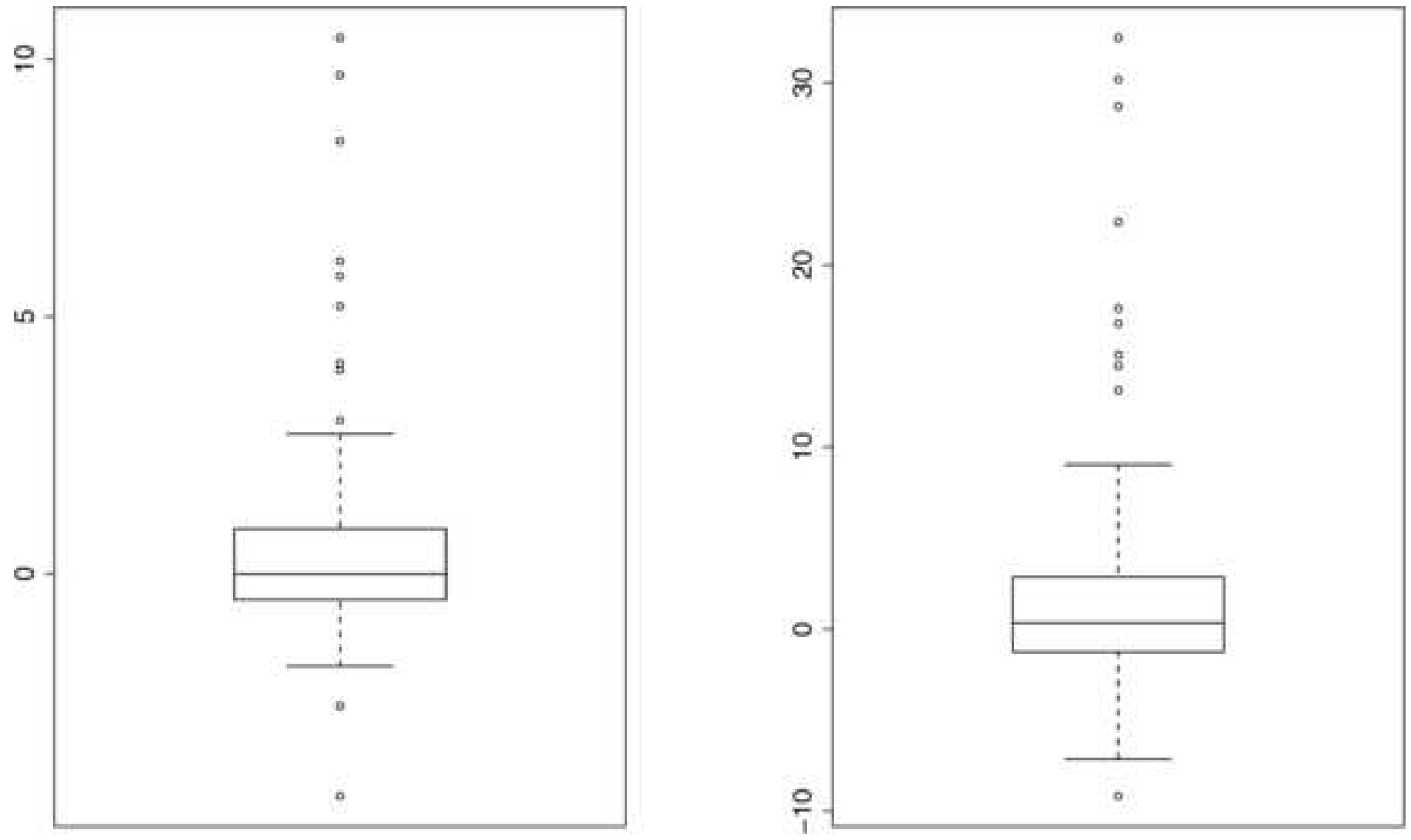

Fig. 8. *Comparison of hard and soft thresholding. Left:* $m(x) = 5x_1^2 x_2^2$, $d = 10$ *and* $\sigma = 0.5$; *right:* $m(x) = 2(x_1 + 1)^3 + 2\sin(10x_2)$, $d = 10$ *and* $\sigma = 1$. *The hard and soft thresholding versions of the rodeo were compared on 100 randomly generated data sets, with a single random test point* $x$ *chosen for each;* $\beta = 0.9$. *The plots show two views of the difference of losses,* $(\widetilde{m}_{\rm hard}(x) - m(x))^2 - (\widetilde{m}_{\rm soft}(x) - m(x))^2$; *positive values indicate an advantage for soft thresholding.*

If $j \in R^c$ then we have $\mathbb{E}(Z_j(x_i)) = o(1)$, so it follows that, conditionally,

$$\mathbb{E}(T_j) = \frac{\sigma^2}{k}\operatorname{tr}(P_j) + o_P(1), \tag{6.8a}$$

$$\operatorname{Var}(T_j) = \frac{2\sigma^4}{k^2}\operatorname{tr}(P_j P_j) + o_P(1). \tag{6.8b}$$

We take the threshold to be

$$\lambda_j = \frac{\widehat{\sigma}^2}{k}\operatorname{tr}(P_j) + 2\frac{\widehat{\sigma}^2}{k}\sqrt{\operatorname{tr}(P_j P_j)\log(n)}. \tag{6.9}$$

Note that if $j > r$, we have

$$\mathbb{E}(T_j) = \frac{1}{k}\sum_i s_j^2(X_i) + O(h_0^6) \tag{6.10}$$

but for $j \le r$ we have

$$\mathbb{E}(T_j) = \frac{1}{k}\sum_i s_j^2(X_i) + O(h_0^2). \tag{6.11}$$

We give an example of this algorithm in the following section, leaving the detailed analysis of the asymptotics of this estimator to future work.

6.4. *Greedier rodeo and LARS.* The rodeo is related to least angle regression (LARS) (Efron et al. [4]). In forward stagewise linear regression, one performs variable selection incrementally. LARS gives a refinement where at each step in the algorithm, one adds the covariate most correlated with the residuals of the current fit, in small, incremental steps. LARS takes steps of a particular size: the smallest step that makes the largest correlation equal to the next-largest correlation. Efron et al. [4] show that the lasso can be obtained by a simple modification of LARS.

The rodeo can be seen as a nonparametric version of forward stagewise regression. Note first that $Z_j$ is essentially the correlation between the $Y_i$'s and the $G_j(X_i, x, h)$s (the change in the effective kernel). Reducing the bandwidth is like adding in more of that variable. Suppose now that we make the following modifications to the rodeo: (i) change the bandwidths one at a time, based on the largest $Z_j^* = Z_j/\lambda_j$, (ii) reduce the bandwidth continuously, rather than in discrete steps, until the largest $Z_j$ is equal to the next largest. Some examples of the greedy version of this algorithm follow.

6.4.1. *Diabetes example.* Figure 9 shows the result of running the greedy version of the rodeo on the diabetes dataset used by [4] to illustrate LARS. The algorithm averages $Z_j^*$ over a randomly chosen set of $k = 100$ data points, and reduces the bandwidth for the variable with the largest value; note that no estimate of $\sigma$ is required. The resulting variable ordering is seen to be very similar to, but different from, the ordering obtained from the parametric LARS fit. The variables were selected in the order 3 (body mass index), 9 (serum), 7 (serum), 4 (blood pressure), 1 (age), 2 (sex), 8 (serum), 5 (serum), 10 (serum), 6 (serum). The LARS algorithm adds variables in the order 3, 9, 4, 7, 2, 10, 5, 8, 6, 1. One notable difference is in the position of the age variable.

6.4.2. *Turlach's example.* In the discussion to the LARS paper, Berwin Turlach [30] gives an interesting example of where LARS and the lasso fails. The function is

$$Y = (X_1 - \tfrac{1}{2})^2 + X_2 + X_3 + X_4 + X_5 + \varepsilon \tag{6.12}$$

with ten variables $X_i \sim \text{Uniform}(0, 1)$ and $\sigma = 0.05$. Although $X_1$ is a relevant variable, it is uncorrelated with $Y$, and LARS and the lasso miss it.

Figure 10 shows the greedy algorithm on this example, where bandwidth corresponding to the largest average $Z_j^*$ is reduced in each step. We use kernel regression rather than local linear regression as the underlying estimator, without first subtracting off a Lasso fit. The variables $x_2$, $x_3$, $x_4$, $x_5$ are linear in the model, but are selected first in every run. Variable $x_1$ is selected fifth in 72 of the 100 runs; a typical run of the algorithm is shown in the left plot. In contrast, as discussed in Turlach [30], LARS selects $x_1$ in position 5 about 25% of the time.

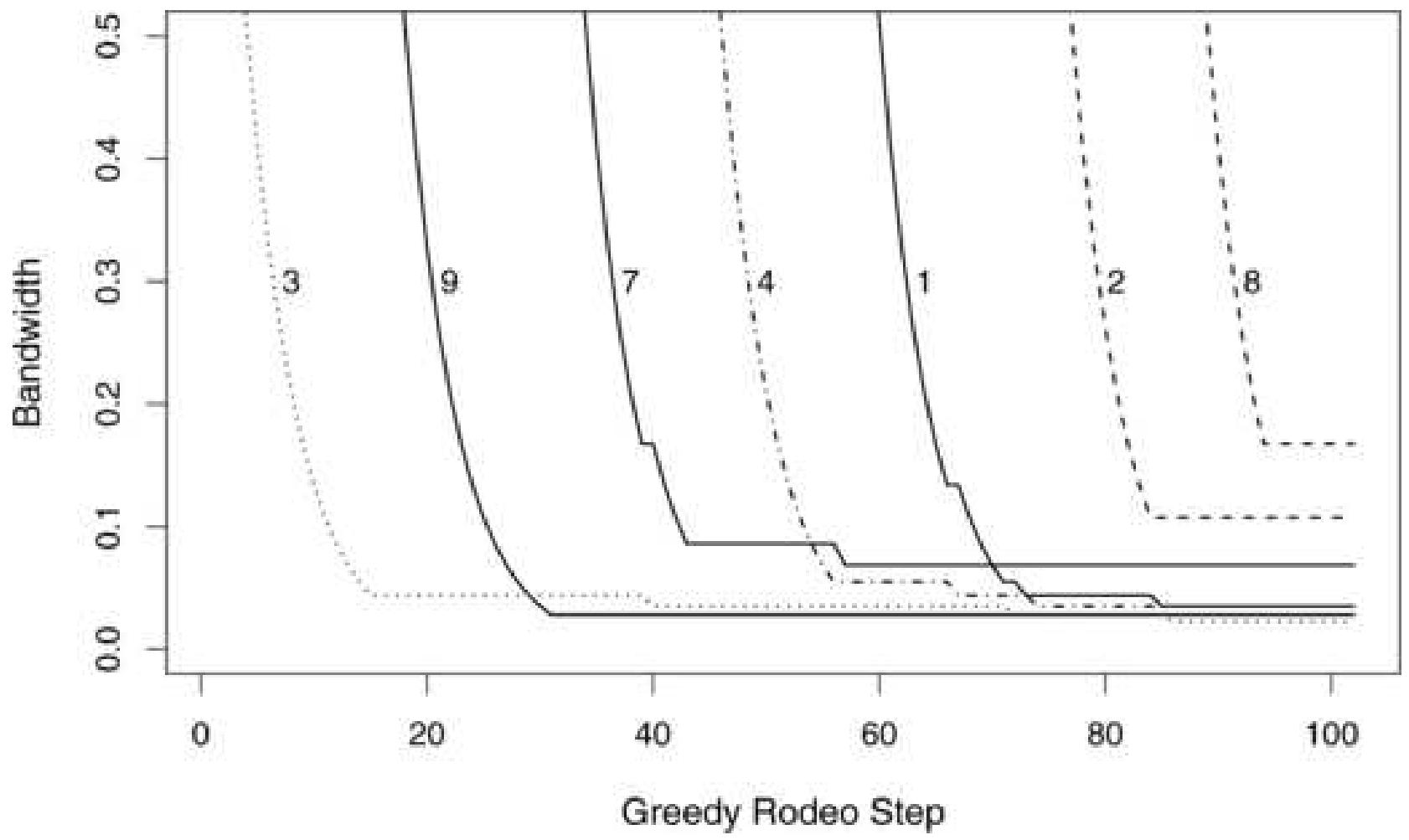

Fig. 9. *Greedy rodeo on the diabetes data, used to illustrate LARS (Efron et al. [4]). A set of $k = 100$ of the total $n = 442$ points were sampled ($d = 10$), and the bandwidth for the variable with largest average $|Z_j|/\lambda_j$ was reduced in each step.*

Figure 11 shows bandwidth traces for this example using the global algorithm described in Section 6.3 with $k = 20$ evaluation points randomly subselected from the data, and $\sigma$ taken to be known. Before starting the rodeo, we subtract off a linear least squares fit, and run the rodeo on the residuals. The first plot shows $h_1, \dots, h_5$. The lowest line is $h_1$ which shrinks the most since $m$ is a nonlinear function of $x_1$. The other curves are the linear effects. The right plot shows the traces for $h_6, \dots, h_{10}$, the bandwidths for the irrelevant variables.

**7. Proofs of technical results.** In this section we give the proofs of the results stated in Section 5. We begin with three lemmas.

We write $Y_n = \widetilde{O}_P(a_n)$ to mean that $Y_n = O_P(b_n a_n)$ where $b_n$ is logarithmic in $n$. As noted earlier, we write $a_n = \Omega(b_n)$ if $\liminf_n |\frac{a_n}{b_n}| > 0$; similarly $a_n = \widetilde{\Omega}(b_n)$ if $a_n = \Omega(b_n c_n)$ where $c_n$ is logarithmic in $n$.

Let

$$\mathcal{H} = \begin{pmatrix} \mathcal{H}_R & 0 \\ 0 & 0 \end{pmatrix}$$

denote the Hessian of $m(x)$. For a given bandwidth $h = (h_1, \dots, h_d)$, denote the bandwidth matrix by $H = \operatorname{diag}(h_1^2, \dots, h_d^2)$. Similarly, let $H_R =$

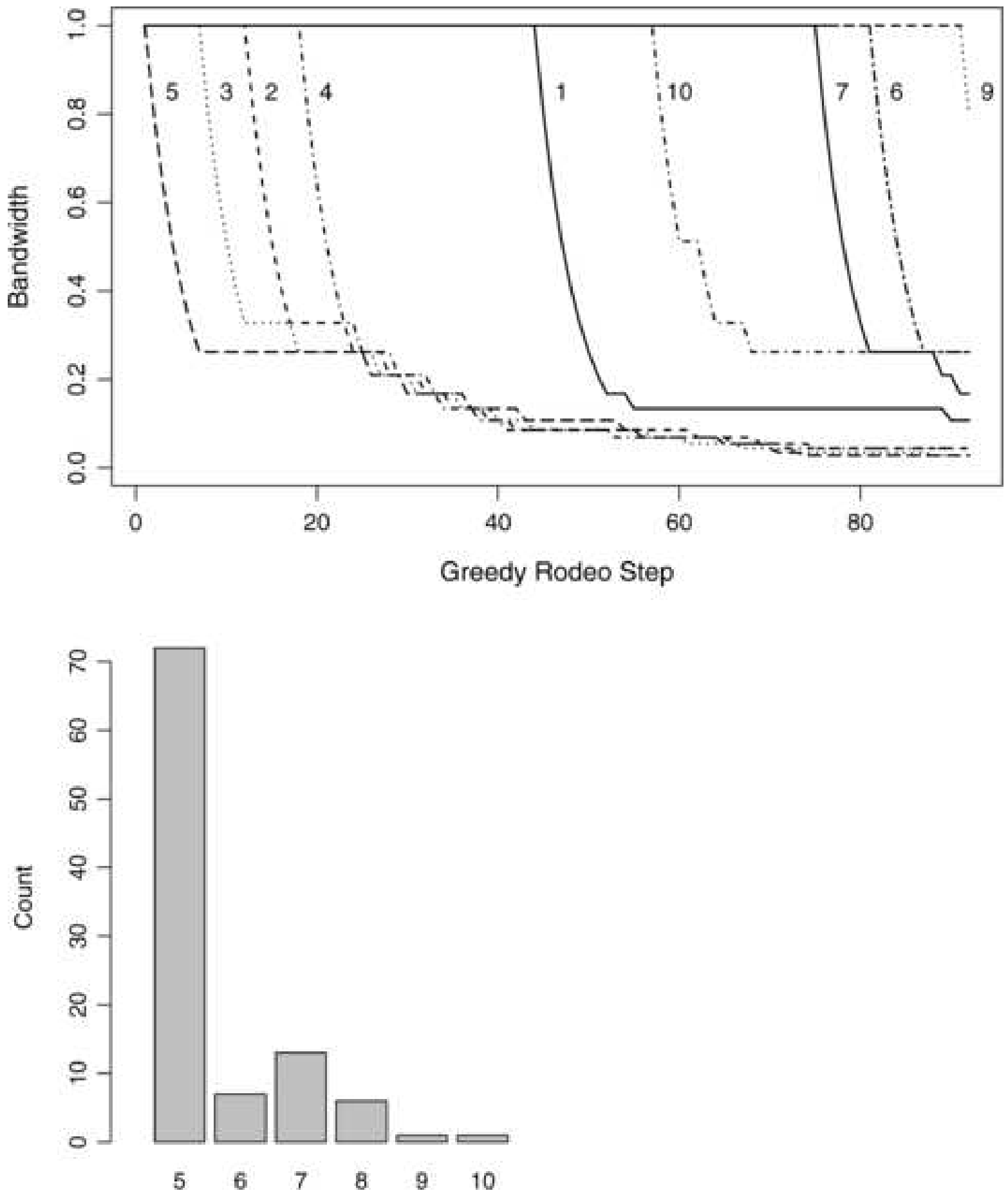

FIG. 10. *Top: A typical run of the greedy algorithm on Turlach's example. The bandwidths are first reduced for variables $x_2, x_3, x_4, x_5$, and then the relevant, but uncorrelated with $Y$ variable $x_1$ is added to the model; the irrelevant variables enter the model last. Bottom: Histogram of the position at which variable $x_1$ is selected, over 100 runs of the algorithm.*

$\mathrm{diag}(h_1^2, \ldots, h_r^2)$. Define

$$\mu_j(h) = \frac{\partial}{\partial h_j} \mathbb{E}[\widehat{m}_H(x) - m(x) | X_1, \ldots, X_n], \tag{7.1}$$

which is the derivative of the conditional bias. The first lemma analyzes $\mu_j(h)$ and $\mathbb{E}(\mu_j(h))$ under the assumption that $f$ is uniform. The second lemma analyzes the variance. The third lemma bounds the probabilities $\mathbb{P}(|Z_j| \geq \lambda_j)$ in terms of tail inequalities for standard normal variables.

In each of these lemmas, we make the following assumptions. We assume that $f$ is uniform, $K$ is a product kernel, and $0 < \beta < 1$. Moreover, we make

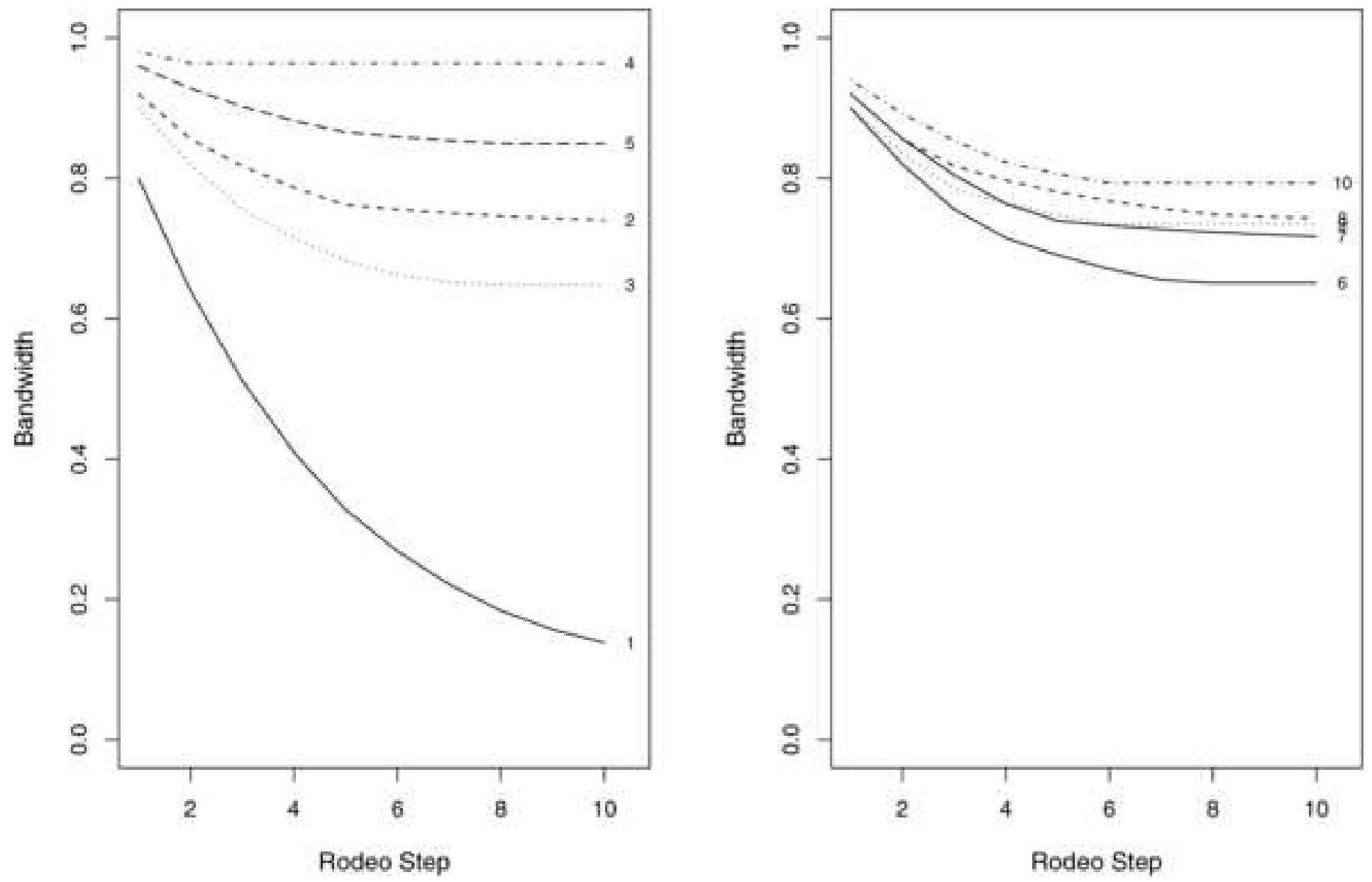

Fig. 11. *The global rodeo averaged over 10 runs on Turlach's example. The left plot shows the bandwidths for the five relevant variables. Since the linear effects (variables two through five) have been subtracted off, bandwidths $h_2, h_3, h_4, h_5$ are not shrunk. The right plot shows the bandwidths for the other, irrelevant, variables.*

use of the following set $\mathcal{B}$ of bandwidths

$$\begin{aligned}\mathcal{B} = \Big\{ h = (h_1, \ldots, h_d) = (\underbrace{\beta^{k_1} h_0, \ldots, \beta^{k_r} h_0}_{r \text{ terms}}, \underbrace{h_0, \ldots, h_0}_{d-r \text{ terms}}) : \\ 0 \le k_j \le T_n,\ j = 1, \ldots, r \Big\},\end{aligned} \tag{7.2}$$

where $T_n \le c_1 \log n$. Finally, we assume that

$$r = O(1), \tag{7.3}$$

$$d = O\left(\frac{\log n}{\log \log n}\right), \tag{7.4}$$

$$h_0 = \frac{c_0}{\log \log n} \qquad \text{for some } c_0 > 0. \tag{7.5}$$

Lemma 7.1. *For each $h \in \mathcal{B}$,*

$$\mathbb{E}(\mu_j(h)) = \begin{cases} \nu_2 m_{jj}(x) h_j + g_j(x_R, h_R) h_j, & j \le r, \\ 0, & j > r, \end{cases} \tag{7.6}$$

*where $\nu_2$ is defined in equation (3.5), and $g_j(x_R, h_R)$ depends only on the relevant variables and bandwidths, and satisfies*

$$|g_j(x_R,h_R)| = O\left(\sum_{k\le r}\sup_x |m_{jjkk}(x)|h_k^2\right). \tag{7.7}$$

*Furthermore, for any $\delta > 0$,*

$$\mathbb{P}\left(\max_{\substack{h\in\mathcal{B}\\1\le j\le d}} \frac{|\mu_j(h)-\mathbb{E}(\mu_j(h))|}{s_j(h)} > \frac{\sqrt{\delta\log n}}{\log\log n}\right) \le \frac{1}{n^{\delta\sigma^2/(8c_0)}} \tag{7.8}$$

*where*

$$s_j^2(h) = \frac{C}{nh_j^2}\prod_{k=1}^d \frac{1}{h_k}, \tag{7.9}$$

*with*

$$C = \sigma^2\int K^2(u)\,du/f(x). \tag{7.10}$$

REMARK 7.2. If we drop the assumption that $f(x) = 1$ then the mean becomes

$$\begin{aligned}&\mathbb{E}(\mu_j(h))\\ &\quad= \begin{cases}\nu_2 m_{jj}(x)h_j + o(h_j), & j\le r,\\ -\operatorname{tr}(H_R\mathcal{H}_R)\nu_2^2(\nabla_j\log f(x))^2h_j + o(h_j\operatorname{tr}(H_R)), & j>r.\end{cases}\end{aligned} \tag{7.11}$$

PROOF OF LEMMA 7.1. We follow the setup of Ruppert and Wand [22] but the calculations need to be uniform over $h\in\mathcal{B}$ and we have to allow for increasing dimension $d$.

Note that there are $N_n = (T_n+1)^r$ elements in the set $\mathcal{B}$. Fix an $h\in\mathcal{B}$. Then $h_j = h_0$ for $j > r$ and

$$\beta^{T_n}h_0 \le h_j \le h_0, \qquad 1\le j\le r. \tag{7.12}$$

Let

$$\mathcal{H} = \begin{pmatrix}\mathcal{H}_R & 0\\ 0 & 0\end{pmatrix}$$

denote the Hessian of $m(x)$. Let $\nabla m$ be the gradient of $m$ at $x$, and let

$$Q = ((X_1-x)^{\mathrm{T}}\mathcal{H}(X_1-x),\ldots,(X_n-x)^{\mathrm{T}}\mathcal{H}(X_n-x))^{\mathrm{T}}. \tag{7.13}$$

Note that $\nabla m$ and $Q$ are only functions of the relevant variables. Then

$$m(X_i) = m(x) + (X_i-x)^{\mathrm{T}}\nabla m + \tfrac{1}{2}Q_i + R_i \tag{7.14}$$

where, using multi-index notation,

$$
\begin{aligned}
R_i &= \tfrac{1}{6}\sum_{|\alpha|=3}(X_i-x)^\alpha \int_0^1 D^\alpha m((1-s)x+sX_i)\,ds \\
&= \sum_{|\alpha|=3}(X_i-x)^\alpha R_\alpha(X_i)
\end{aligned} \tag{7.15}
$$

for functions $R_\alpha$ that only depend on the relevant variables and satisfy

$$|R_\alpha(X_i)| \le \tfrac{1}{6}\sup_x |D^\alpha m(x)|. \tag{7.16}$$

Thus, with $M = (m(X_1), \dots, m(X_n))^{\mathrm{T}}$,

$$M = X_x \begin{pmatrix} m(x) \\ \nabla m \end{pmatrix} + \tfrac{1}{2}Q + R, \tag{7.17}$$

where $R = (R_1, \dots, R_n)^{\mathrm{T}}$. Since $S_x X_x (m(x), \nabla m)^{\mathrm{T}} = m(x)$, the conditional bias

$$b_n(x) = \mathbb{E}(\widehat{m}_H(x)|X_1, \dots, X_n) - m(x) \tag{7.18}$$

is given by

$$b_n(x) = S_x M - m(x) = \frac{1}{2}S_x Q + S_x R \tag{7.19a}$$

$$= \frac{1}{2}e_1^{\mathrm{T}}(X_x^{\mathrm{T}}W_xX_x)^{-1}X_x^{\mathrm{T}}W_xQ + e_1^{\mathrm{T}}(X_x^{\mathrm{T}}W_xX_x)^{-1}X_x^{\mathrm{T}}W_xR \tag{7.19b}$$

$$= \frac{1}{2}e_1^{\mathrm{T}}\Upsilon_n^{-1}\Gamma_n + e_1^{\mathrm{T}}\Upsilon_n^{-1}\frac{1}{n}X_x^{\mathrm{T}}W_xR, \tag{7.19c}$$

where $\Upsilon_n = n^{-1}(X_x^{\mathrm{T}}W_xX_x)$ and $\Gamma_n = n^{-1}(X_x^{\mathrm{T}}W_xQ)$.

*Analysis of* $\Upsilon_n$. We write

$$\Upsilon_n = \begin{pmatrix} \dfrac{1}{n}\sum_{i=1}^n W_i & \dfrac{1}{n}\sum_{i=1}^n W_i(X_i-x)^{\mathrm{T}} \\ \dfrac{1}{n}\sum_{i=1}^n W_i(X_i-x) & \dfrac{1}{n}\sum_{i=1}^n W_i(X_i-x)(X_i-x)^{\mathrm{T}} \end{pmatrix} \tag{7.20a}$$

$$= \begin{pmatrix} A_{11} & A_{12} \\ A_{21} & A_{22} \end{pmatrix}, \tag{7.20b}$$

where

$$W_i = \prod_{j=1}^d \frac{1}{h_j}K\left(\frac{x_j - X_{ij}}{h_j}\right) = \prod_{j=1}^d K_{h_j}(x_j - X_{ij}). \tag{7.21}$$

Now, the variance of $W_i$ can be bounded as

$$\mathrm{Var}(W_i) \le E(W_i^2) = \prod_{j=1}^{d} \frac{1}{h_j} \mathbb{E}\left[\prod_{j=1}^{d} \frac{1}{h_j} K^2\left(\frac{x_j - X_{ij}}{h_j}\right)\right] \tag{7.22}$$

$$= \prod_{j=1}^{d} \frac{1}{h_j} \int K^2(v) f(x + H^{1/2} v)\, dv \tag{7.23}$$

$$= \frac{C}{\sigma^2} \prod_{j=1}^{d} \frac{1}{h_j} \tag{7.24}$$

$$= \frac{n h_j^2}{\sigma^2} s_j^2 \equiv \Delta \tag{7.25}$$

since $f \equiv 1$, where $C$ is as defined in (7.10). Therefore,

$$\mathrm{Var}(A_{11}) \le \frac{\Delta}{n}. \tag{7.26}$$

Also, $W_i \le \Delta$. Hence, by Bernstein's inequality,

$$\mathbb{P}(|A_{11} - \mathbb{E}(A_{11})| > \varepsilon s_j(h)) \le 2 \exp\left\{-\frac{1}{2}\left(\frac{n \varepsilon^2 s_j^2(h)}{\Delta + \Delta \varepsilon s_j(h)/3}\right)\right\} \tag{7.27a}$$

$$= 2 \exp\left\{-\frac{1}{2} \frac{\sigma^2 \varepsilon^2}{h_j^2 (1 + \varepsilon s_j)}\right\}. \tag{7.27b}$$

Now taking $\varepsilon = \frac{\sqrt{\delta \log n}}{\log \log n}$, and using the definition of $h_0 = c_0 / \log \log n$, this gives

$$\mathbb{P}\left(|A_{11} - \mathbb{E}(A_{11})| > s_j(h) \frac{\sqrt{\delta \log n}}{\log \log n}\right) \le 2 \exp\left\{-\frac{\sigma^2 \delta \log n}{4 c_0}\right\} \tag{7.28}$$

$$= 2 n^{-\sigma^2 \delta / (4 c_0)} \tag{7.29}$$

and so with this choice of $\varepsilon$,

$$\begin{aligned} \mathbb{P}\left(\sup_{h \in \mathcal{B}} \frac{|A_{11} - \mathbb{E}(A_{11})|}{s_j(h)} > s_j(h) \frac{\sqrt{\delta \log n}}{\log \log n}\right) &\le 2 (T_n + 1)^r n^{-\sigma^2 \delta / (4 c_0)} \\ &\le n^{-\sigma^2 \delta / (8 c_0)}. \end{aligned} \tag{7.30}$$

Also, since $f(x) = 1$, $\mathbb{E}(A_{11}) = \int \frac{1}{h_1 h_2 \cdots h_d} K(H^{-1/2}(x - u)) f(u)\, du = f(x)$. Hence, for any $\varepsilon > 0$,

$$\mathbb{P}\left(\sup_{h \in \mathcal{B}} \frac{|A_{11} - f(x)|}{s_0(h)} > \varepsilon\right) \to 0. \tag{7.31}$$

Next consider $A_{21}$. Now $\mathbb{E}(A_{21}) = \nu_2(K)HD$ where $D$ is the gradient of $f$. Thus, in the uniform $f$ case, $\mathbb{E}(A_{21}) = \mathbb{E}(A_{12}) = 0$, and by a similar argument as above, $\mathbb{P}(\sup_{h\in\mathcal{B}} |A_{12}|/s_0 > \varepsilon) \to 0$. Turning to $A_{22}$, we again have convergence to its mean and $\mathbb{E}(A_{22}) = \nu_2 f(x)H$.

Thus,

$$\mathbb{E}(\Upsilon_n) = \begin{pmatrix} f(x) & 0 \\ 0 & \nu_2 f(x)H \end{pmatrix}. \tag{7.32}$$

Thus, if

$$\widetilde{\Upsilon}_n^{-1} = \begin{pmatrix} \dfrac{1}{f(x)} & 0 \\ 0 & \dfrac{H^{-1}}{\nu_2 f(x)} \end{pmatrix} \tag{7.33}$$

then

$$\mathbb{P}\left( \max_{jk} \frac{|\Upsilon_n^{-1}(j,k) - \widetilde{\Upsilon}_n^{-1}(j,k))|}{s_0} > \varepsilon \right) \to 0. \tag{7.34}$$

*Analysis of* $\Gamma_n = \frac{1}{n} X_x^{\mathrm{T}} W_x Q$. We can write

$$\begin{aligned} \Gamma_n \equiv \frac{1}{n} X_x^{\mathrm{T}} W_x Q &= \begin{pmatrix} \dfrac{1}{n} \displaystyle\sum_{i=1}^n W_i (X_i - x)^{\mathrm{T}} \mathcal{H}(X_i - x) \\ \dfrac{1}{n} \displaystyle\sum_{i=1}^n (X_i - x) W_i (X_i - x)^{\mathrm{T}} \mathcal{H}(X_i - x) \end{pmatrix} \\ &= \begin{pmatrix} \gamma_1 \\ \gamma_2 \end{pmatrix}. \end{aligned} \tag{7.35}$$

Now,

$$\mathbb{E}(\gamma_1) = \int K(v)(H^{1/2}v)^{\mathrm{T}} \mathcal{H}(H^{1/2}v) f(x + H^{1/2}v)\,dv \tag{7.36a}$$

$$= f(x) \int K(v)(H^{1/2}v)^{\mathrm{T}} \mathcal{H}(H^{1/2}v)\,dv$$

$$+ \int K(v)(H^{1/2}v)^{\mathrm{T}} \mathcal{H}(H^{1/2}v) D^{\mathrm{T}}(H^{1/2}v)\,dv \tag{7.36b}$$

$$+ \tfrac{1}{2} \int K(v)(H^{1/2}v)^{\mathrm{T}} \mathcal{H}(H^{1/2}v)(H^{1/2}v)^{\mathrm{T}} D_2^{\mathrm{T}}(H^{1/2}v)\,dv \tag{7.36c}$$

$$= \nu_2 f(x)\,\mathrm{tr}(H\mathcal{H}_R). \tag{7.36d}$$

The stochastic analysis of $\gamma_1$ from its mean is similar to the analysis of $A_{12}$ and we have $|\gamma_1 - \nu_2 f(x)\,\mathrm{tr}(H\mathcal{H}_R)|/s_0(h) = o(1)$ uniformly.

Next,

$$\mathbb{E}(\gamma_2) = \int K(v)(H^{1/2}v)(H^{1/2}v)^{\mathrm{T}}\mathcal{H}(H^{1/2}v)f(x+H^{1/2}v)\,dv \tag{7.37a}$$

$$\begin{aligned} &= f(x)\int (H^{1/2}v)K(v)(H^{1/2}v)^{\mathrm{T}}\mathcal{H}(H^{1/2}v)\,dv \\ &\quad + \int K(v)(H^{1/2}v)(H^{1/2}v)^{\mathrm{T}}\mathcal{H}(H^{1/2}v)D^{\mathrm{T}}(H^{1/2}v)\,dv \\ &\quad + \tfrac{1}{2}\int K(v)(H^{1/2}v)(H^{1/2}v)^{\mathrm{T}} \\ &\qquad \times \mathcal{H}(H^{1/2}v)(H^{1/2}v)^{\mathrm{T}}D_2^{\mathrm{T}}(H^{1/2}v)\,dv \end{aligned} \tag{7.37b}$$

$$\begin{aligned} &= \int K(v)(H^{1/2}v)(H^{1/2}v)^{\mathrm{T}}\mathcal{H}(H^{1/2}v)D^{\mathrm{T}}(H^{1/2}v)\,dv \\ &\quad + \tfrac{1}{2}\int K(v)(H^{1/2}v)(H^{1/2}v)^{\mathrm{T}} \\ &\qquad \times \mathcal{H}(H^{1/2}v)(H^{1/2}v)^{\mathrm{T}}D_2^{\mathrm{T}}(H^{1/2}v)\,dv \\ &= 0 \end{aligned} \tag{7.37c}$$

Thus,

$$\mathbb{E}(\Gamma_n) = \begin{pmatrix} \nu_2 f(x)\,\mathrm{tr}(H\mathcal{H}_R) \\ 0 \end{pmatrix}. \tag{7.38}$$

*Analysis of remainder* $e_1^{\mathrm{T}}\Upsilon_n^{-1}\frac{1}{n}X_x^{\mathrm{T}}W_xR$. We can write

$$\frac{1}{n}X_x^{\mathrm{T}}W_xR = \begin{pmatrix} \displaystyle\frac{1}{n}\sum_{i=1}^{n}W_iR_i \\ \displaystyle\frac{1}{n}\sum_{i=1}^{n}(X_i-x)W_iR_i \end{pmatrix} = \begin{pmatrix} \delta_1 \\ \delta_2 \end{pmatrix}. \tag{7.39}$$

Then we have, using the definition of $R_\alpha$ in (7.15),

$$\mathbb{E}(\delta_1) = \int K(v)\sum_{|\alpha|=3}R_\alpha(x+H^{1/2}v)(H^{1/2}v)^\alpha f(x+H^{1/2}v)\,dv \tag{7.40}$$

$$= f(x)\sum_{|\alpha|=3}\int K(v)R_\alpha(x+H^{1/2}v)(H^{1/2}v)^\alpha\,dv. \tag{7.41}$$

Due to the fact that $\int K(v)v^\alpha\,dv = 0$ for $\alpha$ of odd order, we can expand one more order to obtain

$$\mathbb{E}(\delta_1) = f(x)\sum_{|\alpha|=4}\int K(v)R_\alpha(x+H^{1/2}v)(H^{1/2}v)^\alpha\,dv \tag{7.42}$$

$$= f(x) \sum_{|\alpha|=4} \bar{R}_\alpha(x,h) h^\alpha \tag{7.43}$$

for functions $\bar{R}_\alpha$ which satisfy

$$|\bar{R}_\alpha(x_R, h_R)| = O\Big(\sup_x |D^\alpha m(x)|\Big). \tag{7.44}$$

Thus, we have that

$$\mathbb{E}(\delta_1) = O\Bigg(\sum_{j,k \le r} \sup_x |m_{jjkk}(x)| h_j^2 h_k^2\Bigg). \tag{7.45}$$

Similarly,

$$\mathbb{E}(\delta_2) = O\Bigg(\sum_{j,j \le r} \sup_x |m_{jkk}(x)| h_j^2 h_k^2\Bigg). \tag{7.46a}$$

Hence,

$$e_1^{\mathrm{T}} \Upsilon_n^{-1} \frac{1}{n} X_x^{\mathrm{T}} W_x R = O_P\Bigg(\sum_{j,k \le r} h_j^2 h_k^2\Bigg). \tag{7.47}$$

Putting all of this together, we conclude that

$$\mathbb{E} b_n(x) = \tfrac{1}{2} \nu_2 \operatorname{tr}(H \mathcal{H}_R) + g(x,h), \tag{7.48}$$

where

$$g(x,h) = O\Bigg(\sum_{j,k \le r} \sup_x |m_{jjkk}(x)| h_j^2 h_k^2\Bigg). \tag{7.49}$$

Taking the derivative with respect to bandwidth $h_j$, for $j \le r$, we obtain

$$\mathbb{E}\mu_j(h) = \nu_2 m_{jj}(x) h_j + g_j(x,h) h_j, \tag{7.50}$$

where

$$g_j(x,h) = O\Bigg(\sum_{k \le r} \sup_x |m_{jjkk}(x)| h_k^2\Bigg). \tag{7.51}$$

The probability bounds established with Bernstein's inequality then give the statement of the lemma. □

REMARK 7.3. Special treatment is needed if $x$ is a boundary point; see Theorem 2.2 of Ruppert and Wand [22].

LEMMA 7.4. *Let* $v_j(h) = \mathrm{Var}(Z_j(h)|X_1,\ldots,X_n)$. *Then*

$$\mathbb{P}\left(\max_{\substack{h\in\mathcal{B}\\1\le j\le d}}\left|\frac{v_j(h)}{s_j^2(h)}-1\right|>\varepsilon\right)\to 0, \tag{7.52}$$

*for all* $\varepsilon > 0$.

PROOF. Let $\ell$ denote the first row of $S_x$. Then, with $\xi \sim N(0,\sigma^2)$,

$$\widehat{m}_H(x) = \sum_i \ell_i Y_i = \sum_i \ell_i m(X_i) + \sum_i \ell_i \varepsilon_i \tag{7.53a}$$

$$\stackrel{d}{=} \sum_i \ell_i m(X_i) + \xi\sqrt{\sum_i \ell_i^2} \tag{7.53b}$$

$$= \sum_i \ell_i m(X_i) + \frac{\Lambda}{\sqrt{nh_1\cdots h_d}}\xi, \tag{7.53c}$$

where

$$\Lambda = \sqrt{nh_1\cdots h_d \sum_i \ell_i^2}. \tag{7.54}$$

Thus,

$$\mathrm{Var}(Z_j(t)|X_1,\ldots,X_n) = \sigma^2\,\mathrm{Var}\left(\frac{\partial}{\partial h_j}\frac{\Lambda}{\sqrt{nh_1\cdots h_d}}\right). \tag{7.55}$$

Now we find an asymptotic approximation for $\Lambda$.

Recall that

$$S_x = \left(\frac{1}{n}X_x^{\mathrm{T}} W_x X_x\right)^{-1}\frac{1}{n}X_x^{\mathrm{T}} W_x \tag{7.56}$$

and from our previous calculations

$$\Upsilon_n^{-1} = \left(\frac{1}{n}X_x^{\mathrm{T}} W_x X_x\right)^{-1} = \begin{pmatrix} \dfrac{1}{f(x)} & 0 \\ 0 & \dfrac{H^{-1}}{\nu_2 f(x)} \end{pmatrix}(1+o_P(1)). \tag{7.57}$$

Note that $\sum_i \ell_i^2$ is the (1, 1) entry of $S_x S_x^{\mathrm{T}}$. But

$$S_x S_x^{\mathrm{T}} = \left(\Upsilon^{-1}\frac{1}{n}X_x^{\mathrm{T}} W_x\right)\left(\Upsilon^{-1}\frac{1}{n}X_x^{\mathrm{T}} W_x\right)^{\mathrm{T}} \tag{7.58a}$$

$$= \frac{1}{n^2}\Upsilon^{-1} X_x^{\mathrm{T}} W_x^2 X_x \Upsilon^{-1} \tag{7.58b}$$

$$= \frac{1}{n}\Upsilon^{-1}\begin{pmatrix} \dfrac{1}{n}\sum_i W_i^2 & \dfrac{1}{n}\sum_i (X_i-x)^{\mathrm{T}} W_i^2 \\ \dfrac{1}{n}\sum_i (X_i-x) W_i^2 & \dfrac{1}{n}\sum_i (X_i-x)(X_i-x)^{\mathrm{T}} W_i^2 \end{pmatrix}\Upsilon^{-1}. \tag{7.58c}$$

So $\Lambda^2$ is the $(1,1)$ entry of

$$\Upsilon^{-1}\begin{pmatrix} \dfrac{h_1\cdots h_d}{n}\displaystyle\sum_i W_i^2 & \dfrac{h_1\cdots h_d}{n}\displaystyle\sum_i (X_i-x)^{\mathrm{T}}W_i^2 \\ \dfrac{h_1\cdots h_d}{n}\displaystyle\sum_i (X_i-x)W_i^2 & \dfrac{h_1\cdots h_d}{n}\displaystyle\sum_i (X_i-x)(X_i-x)^{\mathrm{T}}W_i^2 \end{pmatrix}\Upsilon^{-1}$$

$$= \Upsilon^{-1}\begin{pmatrix} a_{11} & a_{12} \\ a_{21} & a_{22} \end{pmatrix}\Upsilon^{-1}. \tag{7.59}$$

Next, as in our earlier analysis,

$$\mathbb{E}(a_{11}) = \int K^2(v) f(x - H^{1/2}v)\,dv \tag{7.60a}$$

$$= f(x)\int K^2(v)\,dv \tag{7.60b}$$

and similarly, $\mathbb{E}(a_{21}) = \mathbb{E}(a_{22}) = 0$ and $\mathbb{E}(a_{22}) = f(x)\bar{\nu}_2 H$, where $\bar{\nu}_2 I = \int vv^{\mathrm{T}} \times K^2(v)\,dv$. Hence, the leading order expansion of $\Lambda^2$ is given by

$$\frac{\int K^2(v)\,dv}{f(x)} + O(\mathrm{tr}(H)). \tag{7.61}$$

Taking the derivative with respect to $h_j$ we thus conclude that

$$\mathrm{Var}(Z_j(t)|X_1,\ldots,X_n) = \frac{\sigma^2\int K^2(v)\,dv}{f(x)h_j^2}\frac{1}{nh_1\cdots h_d}(1+o_P(1)), \tag{7.62}$$

which gives the statement of the lemma. □

LEMMA 7.5.

1. *For any $c>0$ and each $j>r$,*

$$\mathbb{P}(|Z_j(h_0)| > \lambda_j(h_0)) = o\left(\frac{1}{n^c}\right). \tag{7.63}$$

2. *Uniformly for $h\in\mathcal{B}$ we have the following: for any $c>0$, $j\le r$,*

$$\mathbb{P}(|Z_j(h)| < \lambda_j(h)) \le \mathbb{P}\left(N(0,1) > \frac{\nu_2|m_{jj}(x)|h_j + z_n}{s_j(h)}\right) + o\left(\frac{1}{n^c}\right), \tag{7.64}$$

*where $z_n = O(h_j^3)$.*

PROOF. Proof of (1). Fix $\delta>0, c>0$. By the previous lemmas, there exists a sequence of sets $V_n$ and sequences of constants $\xi_{1,n},\xi_{2,n}$ such that $\xi_{1,n}\le\sqrt{\delta\log n}/\log\log n$, $\xi_{2,n}\to 0$, $\mathbb{P}(V_n^c) = O(n^{-\delta\sigma^2/(8c_0)})$. On $V_n$ we have that

$$|\mu_j(h_0)|/s_j(h_0) \le \xi_{1,n} \tag{7.65}$$

and

$$|s_j(h)/\sqrt{v_j(h_0)}-1| \le \xi_{2,n}. \tag{7.66}$$

The events $V_n$ depend on the $X_i$'s but not $\varepsilon_i$'s. Choosing $\delta$ large enough we have that $\mathbb{P}(V_n^c)=o(n^{-c})$. So, for all $j>r$,

$$\begin{aligned}
&\mathbb{P}(Z_j(h_0)>\lambda_j(h_0)) && \text{(7.67a)}\\
&=\mathbb{P}\biggl(\frac{Z_j(h_0)-\mu_j(h_0)}{s_j(h_0)}>\frac{\lambda_j(h_0)-\mu_j(h_0)}{\sqrt{v_j(h_0)}}\biggr) && \text{(7.67b)}\\
&=\mathbb{E}\biggl(\mathbb{P}\biggl(\frac{Z_j(h_0)-\mu_j(h_0)}{s_j(h_0)}>\frac{\lambda_j(h_0)-\mu_j(h_0)}{\sqrt{v_j(h_0)}}\Big|X_1,\dots,X_n\biggr)\biggr) && \text{(7.67c)}\\
&=\mathbb{E}\biggl(\mathbb{P}\biggl(N(0,1)>\frac{\lambda_j(h_0)-\mu_j(h_0)}{\sqrt{v_j(h_0)}}\Big|X_1,\dots,X_n\biggr)\biggr) && \text{(7.67d)}\\
&=\mathbb{P}\biggl(N(0,1)>\frac{\lambda_j(h_0)-\mu_j(h_0)}{\sqrt{v_j(h_0)}}\biggr) && \text{(7.67e)}\\
&=\mathbb{P}\biggl(N(0,1)>\frac{\lambda_j(h_0)-\mu_j(h_0)}{\sqrt{v_j(h_0)}},V_n\biggr)+o\biggl(\frac{1}{n^c}\biggr) && \text{(7.67f)}\\
&=\mathbb{P}\biggl(N(0,1)>\frac{\lambda_j(h_0)-\mu_j(h_0)}{s_j(h_0)}\frac{\sqrt{v_j(h_0)}}{s_j(h_0)},V_n\biggr)+o\biggl(\frac{1}{n^c}\biggr) && \text{(7.67g)}\\
&=\mathbb{P}(N(0,1)>\sqrt{2\log n}(1-\xi_{2,n})-\xi_{1,n}(1+\xi_{2,n}))+o\biggl(\frac{1}{n^{cl}}\biggr) && \text{(7.67h)}
\end{aligned}$$

and the result follows from the normal tail inequality and the fact that $\sqrt{2\log n}-\xi_{1,n}>\sqrt{(2-\gamma)\log n}$ for any $\gamma>0$.

Proof of (2). By the previous lemmas, there exists a sequence of sets $V_n$ and sequences of constants $\xi_{1,n},\xi_{2,n}$ such that $\xi_{1,n}\le\sqrt{\delta\log n}/\log\log n$, $\xi_{2,n}\to 0$, $\mathbb{P}(V_n^c)=o(1/n^c)$ for any $c>0$ and on $V_n$ we have that $|\mu_j(h)-\nu_2 m_{jj}(x)h_j+O(h_j^3)|/s_j(h)\le\xi_{1,n}$ and $|s_j(h)/\sqrt{v}_j(h)-1|\le\xi_{2,n}$. The events $V_n$ depend on the $X_i$'s but not $\varepsilon_i$'s. Without loss of generality assume that $m_{jj}(x)>0$. Then

$$\begin{aligned}
&\mathbb{P}(|Z_j(h)|<\lambda_j(h)) && \text{(7.68a)}\\
&=\mathbb{P}\biggl(\frac{-\lambda_j(h)-\mu_j(h)}{\sqrt{v_j(h)}}<N(0,1)<\frac{\lambda_j(h)-\mu_j(h)}{\sqrt{v_j(h)}}\biggr) && \text{(7.68b)}\\
&\le\mathbb{P}\biggl(-\infty<N(0,1)<\frac{\lambda_j(h)-\mu_j(h)}{\sqrt{v_j(h)}}\biggr) && \text{(7.68c)}
\end{aligned}$$

$$= \mathbb{P}\left(N(0,1) > \frac{\mu_j(h) - \lambda_j(h)}{\sqrt{v_j(h)}}\right) \tag{7.68d}$$

$$= \mathbb{P}\left(N(0,1) > \frac{\mu_j(h) - \lambda_j(h)}{\sqrt{v_j(h)}} \frac{\sqrt{v}_j(h)}{s_j(h)}\right) \tag{7.68e}$$

$$= \mathbb{P}\left(N(0,1) > \frac{\mu_j(h) - \lambda_j(h)}{\sqrt{v_j(h)}} \frac{\sqrt{v}_j(h)}{s_j(h)}, V_n\right) + o\left(\frac{1}{n^c}\right) \tag{7.68f}$$

$$= \mathbb{P}\left(N(0,1) > \frac{\nu_2 m_{jj}(x) h_j - \lambda_j(h) - z_n}{s_j(h)}(1 - \xi_{2,n})\right) + o\left(\frac{1}{n^c}\right). \quad \square \tag{7.68g}$$

7.1. *Proof of Theorem 5.1.*

Proof of Theorem 5.1. Let $\mathcal{A}_t$ be the active set at step $t$. Define $A_t$ to be the event that $\mathcal{A}_t = \{1, \ldots, r\}$. Let $\mathcal{C}_t = \{\mathcal{A}_t = \varnothing\}$. Recall the definitions of $t_L$ and $t_U$ from equations (5.4) and (5.5). We will show that

$$\mathbb{P}\left(\mathcal{C}_{t_U} \cap \left(\bigcap_{j=1}^{t_L} A_j\right)\right) \to 1 \tag{7.69}$$

from which the theorem follows. We analyze the algorithm as it progresses through steps $1, \ldots, t, \ldots T_n$. Fix $c > 0$. In what follows, we let $\xi_n(c)$ denote a term that is $o(n^{-c})$; we will suppress the dependence on $c$ and simply write $\xi_n$.

*Step* $t = 1$. Define the event

$$B_1 = \{|Z_j| > \lambda_j \text{ for all } j \le r\} \cap \{|Z_j| < \lambda_j \text{ for all } j > r\}. \tag{7.70}$$

Thus, $A_1 = B_1$. We claim that

$$\mathbb{P}(B_1^c) \le \frac{2d}{n} + \xi_n. \tag{7.71}$$

To show (7.71), we proceed as follows. First consider $j > r$. From Lemma 7.5,

$$\mathbb{P}\left(\max_{j>r} |Z_j| > \lambda_j\right) \le \sum_{j=r+1}^{d} \mathbb{P}(|Z_j| > \lambda_j) \le d\xi_n = \xi_n. \tag{7.72}$$

Now consider $j \le r$. Note that $\mu_j^2(h)/s_j^2(h) > 8 \log n$ and hence, from Lemma 7.5,

$$\mathbb{P}(|Z_j| < \lambda_j \text{ for some } j \le r) \le O\left(\frac{1}{n}\right) + \xi_n. \tag{7.73}$$

This proves (7.71).

*Step* $t = 2$. Let $\widetilde{h} = (\widetilde{h}_1, \dots, \widetilde{h}_d)$ be the random bandwidth at step $t = 2$. Let

$$h_* = (\underbrace{\beta h_0, \dots, \beta h_0}_{r \text{ terms}}, \underbrace{h_0, \dots, h_0}_{d-r \text{ terms}}). \tag{7.74}$$

Let $B_2 = \{|Z_1| > \lambda_1, \dots, |Z_r| > \lambda_r\}$. Then $A_2 = B_1 \cap B_2$ and $\widetilde{h} = h_*$ on $A_2$. Now, $\mathbb{P}(A_2^c) \le \mathbb{P}(B_1^c) + \mathbb{P}(B_2^c)$ and

$$\mathbb{P}(B_2^c) = \mathbb{P}(B_2^c, A_1) + \mathbb{P}(B_2^c, A_1^c) \le \mathbb{P}(B_2^c, A_1) + \mathbb{P}(A_1^c) \tag{7.75}$$

$$= \mathbb{P}(B_2^c, A_1) + \frac{1}{n} + \xi_n \tag{7.76}$$

$$= \mathbb{P}\Big(\min_{j \le r} |Z_j(\widetilde{h})| < \lambda_j, A_1\Big) + \frac{1}{n} + \xi_n \tag{7.77}$$

$$= \mathbb{P}\Big(\min_{j \le r} |Z_j(h_*)| < \lambda_j, A_1\Big) + \frac{1}{n} + \xi_n \tag{7.78}$$

$$\le \mathbb{P}\Big(\min_{j \le r} |Z_j(h_*)| < \lambda_j\Big) + \frac{1}{n} + \xi_n \tag{7.79}$$

$$\le \frac{1}{n} + \xi_n + \left(\frac{1}{n} + \xi_n\right) \tag{7.80}$$

where the last step follows from the same argument as in step 1. So,

$$\mathbb{P}(A_2^c) \le \mathbb{P}(A_1^c) + \mathbb{P}(B_2^c) \le 2\mathbb{P}(A_1^c) + \frac{2}{n} + 2\xi_n. \tag{7.81}$$

*Step* $t$ *for* $t \le t_L$. Let $\widetilde{h} = (\widetilde{h}_1, \dots, \widetilde{h}_d)$ be the random bandwidth at step $t$. Let

$$h_* = (\underbrace{\beta^{t-1} h_0, \dots, \beta^{t-1} h_0}_{r \text{ terms}}, \underbrace{h_0, \dots, h_0}_{d-r \text{ terms}}). \tag{7.82}$$

Let $B_t = \{|Z_1| > \lambda_1, \dots, |Z_r| > \lambda_r\}$. Then $A_t = \bigcap_{s=1}^t B_s$ and $\widetilde{h} = h_*$ on $A_t$. Now, $\mathbb{P}(A_t^c) \le \sum_{s=1}^t \mathbb{P}(B_s^c)$ and

$$\mathbb{P}(B_t^c) = \mathbb{P}(B_t^c, A_{t-1}) + \mathbb{P}(B_t^c, A_{t-1}^c) \tag{7.83}$$

$$\le \mathbb{P}(B_t^c, A_{t-1}) + \mathbb{P}(A_{t-1}^c) \tag{7.84}$$

$$= \mathbb{P}(B_2^c, A_{t-1}) + \frac{1}{n} + \xi_n \tag{7.85}$$

$$= \mathbb{P}\Big(\min_{j \le r} |Z_j(\widetilde{h})| < \lambda_j, A_{t-1}\Big) + \frac{1}{n} + \xi_n \tag{7.86}$$

$$= \mathbb{P}\Big(\min_{j \le r} |Z_j(h_*)| < \lambda_j, A_{t-1}\Big) + \frac{1}{n} + \xi_n \tag{7.87}$$

$$\leq \mathbb{P}\left(\min_{j\leq r}|Z_j(h_*)| < \lambda_j\right) + \frac{1}{n} + \xi_n \tag{7.88}$$

$$\leq \frac{1}{n} + \xi_n + \frac{1}{n} + \xi_n \tag{7.89}$$

from Lemma 7.5 and the fact that $\mu_j^2(h)/s_j^2(h) > 8\log n$ for all $t \leq t_U$. Now,

$$\begin{aligned}\mathbb{P}(A_t^c) &\leq \mathbb{P}(A_{t-1}^c) + \mathbb{P}(B_j^c)\\ &\leq \mathbb{P}(A_{t-1}^c) + \left(\mathbb{P}(A_{t-1}^c) + \frac{1}{n}\right) \qquad (7.90)\\ &\leq 2\mathbb{P}(A_{t-1}^c) + \frac{1}{n}\end{aligned}$$

and so, by induction,

$$\begin{aligned}\mathbb{P}(A_t^c) &\leq \frac{2^{t-1}}{n} + \left(\frac{1}{n} + \xi_n\right)\sum_{j=0}^{t-2} 2^j\\ &\leq \frac{2^t}{n} + 2^t\xi_n = o(1)\end{aligned} \tag{7.91}$$

since $2^t\xi_n(c) = o(1)$ for sufficiently large $c$, for all $t \leq t_L$.

*Step* $t = t_U$. Fix $0 < a < 2$. We use the same argument as in the last case except that $\mu_j^2(h)/s_j^2(h) < a\log n$ for $t = t_U$. Let $\chi$ solve $a = 4 - 2\chi - 4\sqrt{1-\chi}$. Then $0 < \chi < 1$ and $\sqrt{2\log n} - \sqrt{a\log n} \geq \sqrt{2(1-\chi)\log n}$. By Lemma 7.5,

$$\begin{aligned}\mathbb{P}(\mathcal{C}_{t_U}) &\leq \mathbb{P}\left(\max_{j\leq r}|Z_j| > \lambda_j\right)\\ &\leq r\mathbb{P}\left(N(0,1) > \sqrt{2(1-\chi)\log n}\right) + \xi_n \qquad (7.92)\\ &\leq \frac{r}{n^{1-\chi}} + \xi_n.\end{aligned}$$

Summarizing,

$$1 - \mathbb{P}\left(\mathcal{C}_{t_U} \cap \left(\bigcap_{s=1}^{t_L} A_s\right)\right) = o(1) \tag{7.93}$$

which proves the theorem. □

PROOF OF COROLLARY 5.2. First note that, for any deterministic bandwidths $h^\star$ satisfying equations (5.2) and (5.3), we have that the squared (conditional) bias is given by

$$\mathrm{Bias}^2(\widehat{m}_{h^\star}) = \left(\sum_{j\leq r} A_j h_j^{\star 2}\right)^2 + o_P(\mathrm{tr}(H^\star)) \tag{7.94a}$$

$$= \sum_{i,j\le r} A_i A_j h_i^{\star 2} h_j^{\star 2} + o_P(\mathrm{tr}(H^\star)) \tag{7.94b}$$

$$= O_P(n^{-4/(4+r)+\varepsilon}) \tag{7.94c}$$

by Theorem 5.1. Similarly, from Theorem 7.4 the (conditional) variance is

$$\mathrm{Var}(\widehat{m}_{h^\star}) = \frac{1}{n}\left(\prod_i \frac{1}{h_i^\star}\right)\frac{R(K)}{f(x)}\sigma^2(1+o_P(1)) \tag{7.95a}$$

$$= O_P(n^{-1+r/(r+4)+\varepsilon}) \tag{7.95b}$$

$$= O_P(n^{-4/(4+r)+\varepsilon}), \tag{7.95c}$$

where $R(K) = \int K(u)^2\,du$. Let $h^\star$ denote the random bandwidths output from the algorithm. There exists sets $V_n$ such that $\mathbb{P}(V_n^c) = o(1)$ and on $V_n$, the bandwidths satisfy equations (5.2) and (5.3). Let $\delta_n = n^{-(4/(4+r))+\varepsilon}$. It follows that

$$\begin{aligned}&\mathbb{P}(|\widehat{m}_{h^\star}(x) - m(x)| > \delta_n)\\ &\quad\le \mathbb{P}(|\widehat{m}_{h^\star}(x) - m(x)| > \delta_n, V_n) + \mathbb{P}(V_n^c) = o(1).\end{aligned} \tag{7.96}$$

□

**Acknowledgments.** We thank Gilles Blanchard and the editors and referees for many helpful comments on this work.

Department of Computer Science
Department of Machine Learning
Carnegie Mellon University
Pittsburgh, Pennsylvania 15213
USA
E-mail: lafferty@cs.cmu.edu

Department of Statistics
Carnegie Mellon University
Pittsburgh, Pennsylvania 15213
USA
E-mail: larry@stat.cmu.edu